\documentclass[a4paper,10pt,leqno,twoside]{article}

\usepackage{amsmath,amsthm,amssymb}
\numberwithin{equation}{section}

\newcommand{\bC}{\mathbb{C}}
\newcommand{\bR}{\mathbb{R}}
\newcommand{\bQ}{\mathbb{Q}}
\newcommand{\bT}{\mathbb{T}}
\newcommand{\bZ}{\mathbb{Z}}

\newcommand{\mf}[1]{\mathfrak{#1}}
\newcommand{\mr}[1]{\mathrm{#1}}
\newcommand{\mcal}[1]{\mathcal{#1}}

\newcommand{\mbb}[1]{\mathbb{#1}}

\def\bm#1{\boldsymbol{#1}} 
\def\per{\operatorname{per}}
\def\tr{\operatorname{tr}}

\def\diag{\operatorname{diag}}
\def\gl{\mathfrak{gl}}

\newcommand{\dete}[1]{\operatorname{det}^{(#1)}}
\newcommand{\imm}[1]{\operatorname{imm}^{#1}}

\def\({ \left( }
\def\){ \right)}

\def\S#1{Section #1}

\theoremstyle{plain}
\newtheorem{thm}{Theorem}[section]
\newtheorem{prop}[thm]{Proposition}
\newtheorem{lem}[thm]{Lemma}
\newtheorem{cor}[thm]{Corollary}
\theoremstyle{definition}

\newtheorem{example}{Example}[section]
\newtheorem{remark}{Remark}[section]
\newtheorem{conj}[thm]{Conjecture}

\def\C{\mathbb{C}}
\newcommand{\U}[1][n]{\mathcal{U}(\gl_{#1})} 
\def\Znn{\mathbb{Z}_+} 
\newcommand{\GLmod}[2][n]{\mathcal{M}_{#1}^{#2}}
\newcommand{\M}[2]{\mathbb{M}_{#1,#2}}
\newcommand{\tM}[2]{\widetilde{\mathbb{M}}_{#1,#2}}
\def\st{\mathcal{S}} 
\def\kakko#1{\left(#1\right)}
\def\ckakko#1{\left\{#1\right\}}
\def\floor#1{\left\lfloor#1\right\rfloor}
\def\then{\,\Longrightarrow\,}
\newenvironment{smatrix}
{\left(\begin{smallmatrix}}{\end{smallmatrix}\right)}
\def\eqspace{\phantom{{}={}}}
\newcommand{\adet}[1][\alpha]{\operatorname{det}^{(#1)}}
\def\sym#1{\mathfrak{S}_{#1}}
\newcommand{\A}[1][n]{\mathcal{A}(\operatorname{Mat}_{#1})}
\newcommand{\rpn}[1][n]{\rho\kern-0.15em\lower0.8ex\hbox{$\scriptstyle\gl_{#1}$}\!}
\def\len#1{\ell(#1)}
\def\rk{\operatorname{rk}}
\date{\today}
\newenvironment{keywords}{\smallskip\noindent{\bfseries Keywords:}}{}
\newenvironment{MSC}{\smallskip\noindent{\bfseries 2000 Mathematical Subject Classification:}}{}
\def\Uq{\mathcal{U}_q(\mathfrak{gl}_n)}

\DeclareMathOperator{\ind}{ind}

\DeclareMathOperator{\End}{End}
\DeclareMathOperator{\Mat}{Mat}
\DeclareMathOperator{\rank}{rk}
\DeclareMathOperator{\fix}{fix}

\def\reftab{\mathbb{T}}

\def\N{\mathbb{N}}
\def\C{\mathbb{C}}

\def\sym#1{\mathfrak{S}_{#1}}

\def\card#1{\left|#1\right|}
\def\kakko#1{\left(#1\right)}
\def\ckakko#1{\left\{#1\right\}}

\def\ceil#1{\left\lceil#1\right\rceil}

\newcommand{\inprod}[3][{}]{\left\langle#2,\,#3\right\rangle_{\!#1}}

\def\st{\operatorname{Sym}}

\def\len#1{\ell(#1)}
\def\deq{=}

\def\triv#1{\boldsymbol{1}_{#1}}

\def\bc{\overline{c}}

\def\phi{\varphi}

\newcommand{\tramat}[2][n,l]{F^{#2}_{#1}}
\newcommand{\gcp}[3][\alpha]{f^{#2}_{#3}(#1)}

\newcommand{\cp}[1]{f_{#1}}

\def\ys#1{c_{\lower0.3ex\hbox{$\scriptstyle#1$}}}
\def\bys#1{\bc_{\lower0.3ex\hbox{$\scriptstyle#1$}}}

\newcommand{\cw}[1]{\nu(#1)} 
\newcommand{\anu}[1]{\alpha^{\cw{#1}}}
\newcommand{\ve}{\boldsymbol{e}}

\newcommand{\Smod}[1]{\mathcal{S}^{#1}}
\newcommand{\Cmod}[2]{\boldsymbol{V}_{\!\!#1,#2}}

\newcommand{\sfphi}[2][K]{\phi^\sharp_{\!\lower0.4ex\hbox{$\scriptstyle#1,#2$}}}

\newcommand{\intertwiner}{\mathcal{T}}

\newcommand{\tHGF}[5]
{{}_{#1}\tilde{F}_{#2}\kakko{\genfrac{}{}{0pt}{}{#3}{#4};#5}}
\def\phs#1#2{(#1)_{#2}}

\theoremstyle{definition}

\newtheorem{ex}[thm]{Example}
\newtheorem*{ackn}{Acknowledgement}

\theoremstyle{remark}
\newtheorem{rem}[thm]{Remark}

\newenvironment{dedication}{\begin{center}\itshape}{\end{center}\medskip}

\title{\bfseries
Alpha-determinant cyclic modules \\ and Jacobi polynomials}
\author{\scshape Kazufumi KIMOTO, %
Sho MATSUMOTO\thanks{%
Research Fellow of the Japan Society 
for the Promotion of Science, partially supported by Grant-in-Aid 
for Scientific Research (C) No. 17006193.}
\,and 
Masato WAKAYAMA\thanks{%
Partially supported by Grant-in-Aid
for Exploratory Research No. 18654005.},\\
\scshape and with an Appendix by \textsc{Kazufumi KIMOTO}
}
\date{May 29, 2008}
\begin{document}

\maketitle

\begin{dedication}
Dedicated to Professor Masaaki Yoshida on his sixtieth birthday
\end{dedication}

\begin{abstract}
For positive integers $n$ and $l$,
we study the cyclic $\U$-module generated by the $l$-th power of the
$\alpha$-determinant $\dete{\alpha}(X)$.
This cyclic module is isomorphic to
the $n$-th tensor space $(\st^l(\C^n))^{\otimes n}$
of the symmetric $l$-th tensor space of $\bC^n$
for all but finite exceptional values of $\alpha$.
If $\alpha$ is exceptional,
then the cyclic module is equivalent to
a \emph{proper} submodule of $(\st^l(\C^n))^{\otimes n}$,
i.e. the multiplicities of several irreducible subrepresentations
in the cyclic module are smaller
than those in $(\st^l(\C^n))^{\otimes n}$.
The degeneration of each isotypic component of the cyclic module
is described by a matrix
whose size is given by a Kostka number
and entries are polynomials in $\alpha$ with rational coefficients.
Especially, we determine the matrix completely when $n=2$.
In that case, the matrix becomes a scalar
and is essentially given by the classical Jacobi polynomial.
Moreover, we prove that these polynomials are unitary.

In the Appendix, we consider 
a variation of the spherical Fourier transformation for $(\sym{nl},\sym l^n)$ as a main tool
to analyze the same problems,
and describe the case where $n=2$ 
by using the zonal spherical functions of the Gelfand pair $(\sym{2l},\sym l^2)$.

\begin{keywords}
Alpha-determinant, cyclic modules, Jacobi polynomials, 
singly confluent Heun ODE, 
permanent, Kostka numbers, irreducible decomposition,
spherical Fourier transformation,
zonal spherical functions, Gelfand pair.
\end{keywords}

\begin{MSC}
22E47, 
33C45, 
43A90. 
\end{MSC}
\end{abstract}


%
\section{Introduction}
%

Let $\A$ be the associative $\C$-algebra consisting
of polynomials in variables $\{x_{ij}\}_{1 \le i,j \le n}$.
We introduce a $\U$-module structure on $\A$,
where $\U$ is the universal enveloping algebra of
the general linear Lie algebra $\gl_n=\gl_n(\C)$, by
$$
\rpn(E_{ij})f=\sum_{k=1}^n x_{ik} \frac{\partial f}{\partial x_{jk}}
\qquad\kakko{f \in \mcal{A}(\mathrm{Mat}_n)},
$$
which is obtained as a differential representation of the translation of $GL_n=GL_n(\C)$.

Since the \emph{determinant} $\det(X)$ of the matrix
$X=(x_{ij})_{1\le i,j\le n}$ is a relative $GL_n$-invariant in $\A$,
obviously the linear span $\C\cdot\det(X)$ is
a one-dimensional irreducible (highest weight) $\U$-submodule of $\A$.
This submodule is equivalent to
the skew-symmetric tensor representation $\wedge^n(\C^n)$ of
the natural representation of $\U$ on $\C^n$.
The symmetric counterpart of the determinant is
the \emph{permanent} $\per(X)$ given by
\begin{equation*}
\per(X)=\sum_{\sigma\in\sym n}x_{\sigma(1)1}x_{\sigma(2)2}\cdots x_{\sigma(n) n}.
\end{equation*}
Although $\per(X)$ is not a relative invariant of $GL_n$,
the cyclic module $\rpn(\U)\cdot\per(X)$
(i.e. the smallest invariant subspace in $\A$ containing $\per(X)$)
is irreducible and is equivalent to the symmetric tensor representation $\st^n(\C^n)$
of the natural representation.

The \emph{$\alpha$-determinant} of $X$ is defined by
\begin{equation}\label{eq:alpha_det}
\dete{\alpha}(X)=\sum_{\sigma \in \mf{S}_n}\alpha^{\nu(\sigma)}
x_{\sigma(1)1}x_{\sigma(2)2}\cdots x_{\sigma(n) n},
\end{equation}
where $\nu(\sigma)$ is
\emph{$n$ minus the number of cycles in $\sigma \in \mf{S}_n$}.
The notion of the $\alpha$-determinant was first introduced in \cite{VereJones}
in order to describe the coefficients in the expansion of $\det(I-\alpha A)^{-1/\alpha}$,
which is used to treat
the multivariate binomial and negative binomial distributions
in a unified way.
Later, it is also used to define a certain random process in \cite{ST}.
We note that a pfaffian analogue ($\alpha$-pfaffian) is also introduced
and studied in the same (probability theoretic) view point
by the second author in \cite{Mat}.

The $\alpha$-determinant is a common generalization of
(and/or an interpolation between)
the determinant and permanent since
$\adet[-1](X)=\det(X)$ and $\adet[1](X)=\per(X)$.
In this sense, the $\alpha$-determinant cyclic module
$\rpn(\U)\cdot\adet(X)$ is regarded as an interpolation
of two irreducible representations --- the skew-symmetric tensor representation
and symmetric tensor representation.
In \cite{MW}, the second and third authors determined
the structure of the $\U$-cyclic module $\rpn(\U)\cdot\adet(X)$.
The irreducible decomposition of $\rpn(\U)\cdot\adet(X)$ is given by
\begin{equation}\label{eq:result_of_MW}
\rpn(\U)\cdot\adet(X)\cong%
\bigoplus_{\substack{\lambda\vdash n\\ f_\lambda(\alpha)\ne0}}
\kakko{\GLmod\lambda}^{\oplus f^\lambda}.
\end{equation}
Here we denote by $\GLmod\lambda$ the irreducible highest weight $\U$-module
of highest weight $\lambda$
(we identify the highest weight and the corresponding partition),
$f^\lambda$ the number of standard tableaux with shape $\lambda$
and $f_\lambda(\alpha)$ the (modified) content polynomial \cite{Mac}
for $\lambda$ defined by
\begin{equation}\label{Polynomial_f}
f_\lambda(\alpha)=\prod_{i=1}^{\len\lambda}\prod_{j=1}^{\lambda_i}(1+(j-i)\alpha).
\end{equation}
In other words,
the structure of $\rpn(\U)\cdot\adet(X)$ changes drastically
when $\alpha=\pm1/k$ ($k=1,2,\dots,n-1$).
This result implies that
$\adet(X)$ may obtain some special feature like
$\det(X)$ and/or $\per(X)$ for such special values of $\alpha$.
Actually,
when $\alpha=-1/k$ for some $k$,
$\adet[-1/k](X)$ has an analogous property
of the alternating property of the determinant.
Based on this fact, for instance,
we can construct a relative $GL_n$-invariant from $\adet[-1/k](X)$
(see \cite{KWinvariant}).
It is worth noting that
we also introduced
an analogous object of the $\alpha$-determinant
$\adet_q(X)$
in the quantum matrix algebra,
and study the quantum enveloping algebra cyclic module
$\Uq\cdot\adet_q(X)$
in \cite{KWquantum}.
Compared to the classical case \cite{MW},
the cyclic module in the quantum case is much complicated whereas
has a rich structure.

As a next stage,
as in the beginning of the study
of infinite dimensional representation theory by Gel'fand 
and Na{\u\i}mark \cite{GN} in the middle of the last century,
it is natural to proceed in the study of
the cyclic modules $\rpn(\U)\cdot\adet(X)^s$ for $s\in\C$
under a suitable reformulation (see \S\ref{Subsection:ComplexCase}).
In this case,
the cyclic modules $\rpn(\U)\cdot\adet(X)^s$
is not finite dimensional in general.
Actually,
if $s$ is \emph{not} a nonnegative integer,
then $\adet(X)^s$ is no longer a polynomial
and $\rpn(\U)\cdot\adet(X)^s$ becomes infinite dimensional
unless $\alpha=-1$.
On the contrary,
when $s=l$ is a \emph{positive integer},
$\rpn(\U)\cdot\adet(X)^l$ is a submodule of the polynomial algebra $\A$
and is \emph{finite dimensional}.

In this article,
we treat the finite-dimensional cases,
that is, we study the cyclic module
$\rpn(\U)\cdot\adet(X)^l$ for a given positive integer $l$.
We first show that the irreducible decomposition is given in the form
\begin{equation}
\rpn(\U)\cdot\dete{\alpha}(X)^l
\cong \bigoplus_{\lambda \vdash nl} 
(\mcal{M}^{\lambda}_n)^{\oplus m^\lambda_{n,l}(\alpha)},
\end{equation}
where $m^\lambda_{n,l}(\alpha)$ denotes the multiplicity
of the irreducible submodule with highest weight $\lambda$
which satisfies $0 \le m^\lambda_{n,l}(\alpha) \le K_{\lambda (l^n)}$
(Theorem \ref{thm:decomposition}).
Here $K_{\lambda \mu}$ is the \emph{Kostka number}
defined as the number of semi-standard tableaux
of shape $\lambda$ and weight $\mu$.
Moreover,
there exists a certain matrix $\tramat\lambda(\alpha)$
of size $K_{\lambda (l^n)}$,
which is called the \emph{transition matrix} for $\lambda$,
whose entries are polynomials in $\alpha$
such that $m^\lambda_{n,l}(\alpha)=\rk\tramat\lambda(\alpha)$
for each $\lambda$.
By this fact,
for all but finitely many $\alpha$,
we have $m^\lambda_{n,l}(\alpha)=K_{\lambda(l^n)}$ for any $\lambda$.
Namely, the cyclic module
$\rpn(\U)\cdot\adet(X)^l$
is equivalent to the space
$(\st^l(\C^n))^{\otimes n}$
of the symmetric $l$-tensors on $\bC^n$ for almost all $\alpha$.
We note that
$\tramat[n,1]\lambda(\alpha)$ is a scalar matrix
$f^\lambda(\alpha)\cdot I$
(see \eqref{eq:result_of_MW} and \eqref{Polynomial_f}).

Consequently,
we have to describe
the transition matrix $\tramat\lambda(\alpha)$
and/or its rank
$\rk\tramat\lambda(\alpha)(=m^\lambda_{n,l}(\alpha))$
explicitly.
When $n=2$,
we can completely determine the explicit form
of the transition matrices (see Section 4).
In this case,
each transition matrix 
is a scalar and given by a classical \emph{Jacobi polynomial}.
(Precisely, the scalar satisfies a \emph{singly confluent Heun ordinary differential equation} with respect to $\alpha$.
See Corollary \ref{cor:Heun}.)
In other words,
the Jacobi polynomials play the role of the content polynomials.
Moreover, one shows that
these Jacobi polynomials are \emph{unitary},
and hence the multiplicity $m_{2,l}^\lambda(\alpha)$
is non-zero unless $|\alpha|=1$
for each partition $\lambda$ of $2l$.
These are our main result.

Here we should remark that
the Jacobi polynomial does \emph{not} appear
as a spherical function (i.e. a matrix coefficient of a representation) 
in our story,
and hence, it is important to clarify the reason
why the transition matrix becomes a (unitary) Jacobi polynomial when $n=2$.
It seems a far-reaching matter at present
to describe the transition matrices when $n\ge3$.
In fact, we can only give
explicit expressions of transition matrices
in a few special cases.
It is not clear whether
(the entries of) the transition matrices are given by
certain special polynomials.
We leave these problems to the future study.

This paper is organized as follows. 
In \S 2, we recall the $GL_n$-module structure of the 
tensor space $(\st^l(\bC^n))^{\otimes n}$.
This space is the basic one for the study of 
$\alpha$-determinant cyclic modules.
In \S 3, we study the structure of the cyclic module $\U \cdot \dete{\alpha}(X)^l$.
The transition matrix, which determines the multiplicity of the irreducible component
in the cyclic module, is defined in this section.
In \S 4, we exclusively deal with the simple case where $n=2$. 
As stated above, the transition matrix in this case is explicitly given 
by a classical Jacobi polynomial.
In \S 5, we give a conjecture for the
permanent cyclic module $\U\cdot\per(X)^l$ ($\alpha=1$ case),
introduce a certain suitable reformulation of our problem
for the general complex power cases (i.e. $\U\cdot\adet(X)^s$ for $s\in\C$)
and give a remark on $\phi$-immanant cyclic modules
which is a generalization of the situation.

In the Appendix, we investigate our problems by another approach;
We adopt a variation of the spherical Fourier transformation for $(\sym{nl},\sym l^n)$ as a main tool
to analyze the structure of $\U \cdot \dete{\alpha}(X)^l$,
and give another proof of the results in \S 3.
We also describe the transition matrices in the case where $n=2$ 
by using the zonal spherical functions of the Gelfand pair $(\sym{2l},\sym l^2)$.


%
\section{Preliminaries on representation of $\U$}
%

Let $\Znn$ be the set of all non-negative integers.
For a positive integer $n$, we put $[n]=\{1,2,\dots, n\}$.
Let $\bm{e}_1, \dots, \bm{e}_n$ be the standard basis of $\bC^n$. 
The symmetric $l$-tensor space $\mathcal{S}^l(\bC^n)$ is 
the set of all polynomials of degree $l$
in variables $\bm{e}_i$ and expressed as follows:
$$
\st^l(\bC^n) = \bigoplus_{\begin{subarray}{c} m_1, \dots, m_n \in \bZ_+, \\
m_1+ \cdots +m_n =l \end{subarray}} \bC \cdot 
\bm{e}_{1}^{m_1} \bm{e}_{2}^{m_2} \cdots \bm{e}_{n}^{m_n}.
$$

Let $\M nl$ be the set of all $\Znn$-matrices of size $n$
such that 
the sum of entries in each column is equal to $l$:
$$
\M nl=\Bigl\{M=(m_{ij})_{1 \le i,j \le n} \ \Big| \ m_{ij} \in \bZ_+, \ 
\sum_{i=1}^n m_{ij}=l \ (1 \le j \le n)\Bigr\}.
$$
Put 
$$
\bm{e}^M = \bm{e}_{1}^{m_{11}} \bm{e}_{2}^{m_{21}} \cdots \bm{e}_{n}^{m_{n1}} \otimes 
\cdots \otimes  
\bm{e}_{1}^{m_{1n}} \bm{e}_{2}^{m_{2n}} \cdots \bm{e}_{n}^{m_{nn}}
$$
for each $M \in \M nl$.
Then the tensor space $(\st^l (\bC^n))^{\otimes n}$ is given by
$$
(\st^l (\bC^n))^{\otimes n} = \bigoplus_{M \in \M nl} \bC \cdot \bm{e}^M.
$$

The universal enveloping algebra $\U$ acts on $\bC^n$ in a natural way:
$E_{ij} \cdot \bm{e}_k = \delta_{jk} \bm{e}_i$,
where $\delta_{jk}$ is Kronecker's delta.
This action induces the action of
$\U$ on $(\st^l(\bC^n))^{\otimes n}$ as
\begin{equation} \label{eq:ActionSymTensor}
E_{pq} \cdot \bm{e}^M = \sum_{k=1}^n m_{qk} \bm{e}^{M+R^{pq}_k}
\qquad (1 \le p,q \le n, \ 
M=(m_{ij})_{1 \le i,j \le n} \in \mathbb{M}_{n,l}),
\end{equation}
where $R^{pq}_k$ is the matrix of size $n$ whose $(i,j)$-entry is equal to 
$(\delta_{ip} - \delta_{iq}) \delta_{jk}$.
We note that $R_k^{pq}=-R_k^{qp}$.
The irreducible decomposition of the 
$\U$-module $(\st^l(\bC^n))^{\otimes n}$ 
is well known and given by
$$
(\st^l(\bC^n))^{\otimes n} \cong \bigoplus_{\lambda \vdash nl}
 (\GLmod\lambda)^{\oplus K_{\lambda (l^n)}},
$$
see e.g. \cite{FultonHarris, Weyl}.
Here $\GLmod\lambda$ denotes the highest weight module of $\U$ with 
highest weight $\lambda=(\lambda_1,\dots, \lambda_n)$ and
$K_{\lambda (l^n)}$ denotes the Kostka number
which is defined as
the number of semi-standard tableaux of shape $\lambda$ and weight 
$(l^n)= (l,l,\dots,l)$.

\begin{example}
Let $n=2$. Then 
$\mathbb{M}_{2,l} = \left\{\begin{smatrix} r & s \\ l-r & l-s \end{smatrix}
\ \big| \ 0 \le r,s \le l 
\right\}$.
When $l=2$ we see that
\begin{align*}
E_{21} \cdot \bm{e}^{\begin{smatrix} 2 & 1 \\ 0& 1 \end{smatrix}} =& 2
\bm{e}^{\begin{smatrix} 2 & 1 \\ 0& 1 \end{smatrix}+ 
\begin{smatrix} -1 & 0 \\ 1& 0 \end{smatrix}} +
\bm{e}^{\begin{smatrix} 2 & 1 \\ 0& 1 \end{smatrix}+ 
\begin{smatrix} 0 & -1 \\ 0& 1 \end{smatrix}} = 2
\bm{e}^{\begin{smatrix} 1 & 1 \\ 1& 1 \end{smatrix}}
+\bm{e}^{\begin{smatrix} 2 & 0 \\ 0& 2 \end{smatrix}}
=2 \bm{e}_1 \bm{e}_2 \otimes \bm{e}_1 \bm{e}_2 + \bm{e}_1^2 \otimes \bm{e}_2^2, \\
E_{11} \cdot \bm{e}^{\begin{smatrix} 2 & 1 \\ 0& 1 \end{smatrix}} =& 
 3
\bm{e}^{\begin{smatrix} 2 & 1 \\ 0& 1 \end{smatrix}}
=3 \bm{e}_1^2 \otimes \bm{e}_1 \bm{e}_2
\end{align*}
for instance.
The irreducible decomposition of $(\st^l(\bC^2))^{\otimes 2}$ is given as
$$
(\st^l(\bC^2))^{\otimes 2} \cong
\bigoplus_{s=0}^l \GLmod[2]{(2l-s,s)}.
$$
\qed
\end{example}

The following lemma plays a fundamental role
in the discussion below.

\begin{lem} \label{lem:Cyclic}
Let $I_n$ be the identity matrix of size $n$.
Then it holds that
$(\st^l(\bC^n))^{\otimes n}=\U\cdot\bm{e}^{lI_n}$.
Namely, the vector $\bm{e}^{l I_n} 
=\bm{e}_1^{l} \otimes \bm{e}_2^l \otimes \cdots \otimes \bm{e}_n^l$ 
is a cyclic vector
of the $\U$-module $(\st^l(\bC^n))^{\otimes n}$.
\end{lem}

\begin{proof}
Fix a positive integer $l$.
Let $\tM nl$ be the subset
$$
\tM nl = \{ M \in \M nl \ | \ \bm{e}^M\in\U\cdot\bm{e}^{l I_n} \}
$$
of $\M nl$.
Let us prove $(\st^l(\bC^n))^{\otimes n}= 
\bigoplus_{M \in \tM nl} \bC \cdot \bm{e}^M$,
or equivalently
\begin{equation} \label{eq:MathbbM}
\widetilde{\mathbb{M}}_{n,l} \supset \mathbb{M}_{n,l}
\end{equation}
by induction on $n$. 

The universal enveloping algebra $\U[n-1]$ is
embedded in $\U$ as a subalgebra
in a natural way. 
Assume that the inclusion \eqref{eq:MathbbM} holds up to $n-1$.
Then the matrices of the form
$$
M' \oplus (l)=
\( \begin{array}{@{\,}ccc|c@{\,}}
 & & & 0 \\
 & M' & & \vdots \\
 & & & 0 \\ \hline
0 & \ldots & 0 & l \end{array} \)
\qquad (\text{$M' \in \mathbb{M}_{n-1,l}$})
$$
are contained in $\widetilde{\mathbb{M}}_{n,l}$
by the introduction hypothesis.
Applying several $E_{jn}$'s successively
to $\bm{e}^{M'\oplus (l)}$ suitably many times,
we first see that 
any matrix of the form 
\begin{equation} \label{eq:MatrixInduction}
\( \begin{array}{@{\,}ccc|c@{\,}}
 & & & m_{1n} \\
 & M' & & \vdots \\
 & & & m_{n-1,n} \\ \hline
0 & \ldots & 0 & m_{nn} \end{array} \) \qquad 
(M' \in \mathbb{M}_{n-1,l}, \ \sum_{i=1}^n m_{in}=l)
\end{equation}
belongs to $\widetilde{\mathbb{M}}_{n,l}$.

Next, we put 
$$
\M nl(p) = \Bigl\{M \in \mathbb{M}_{n,l} \ \Big| \ 
j \in \{p+1, p+2,\dots, n-1\} \text{\,\ and\,\ } i \in \{1,2,\dots, j-1, n\}
\then m_{ij}=0 \Bigr\}
$$
for each $0 \le p \le n-1$.
Notice that
$\M nl(0)\subset\M nl(1)\subset\cdots\subset\M nl(n-1)=\M nl$.
We show that $\mathbb{M}_{n,l}(p) \subset \widetilde{\mathbb{M}}_{n,l}$
for any $0 \le p \le n-1$ by induction on $p$.
By definition, we see that
any element in $\mathbb{M}_{n,l}(0)$ is of the form \eqref{eq:MatrixInduction},
so that we have $\mathbb{M}_{n,l}(0) \subset \widetilde{\mathbb{M}}_{n,l}$. 
Assume $\mathbb{M}_{n,l}(p-1) \subset \widetilde{\mathbb{M}}_{n,l}$ for $1\le p\le n-1$.
Take any matrix $M=(m_{ij})_{1 \le i,j \le n}$ in $\mathbb{M}_{n,l}(p)$,
and put $\widetilde{M}=M+ \sum_{i=1}^{p-1} m_{ip} R_{p}^{pi} + m_{np} R_{p}^{pn}$.
Equivalently, 
$\widetilde{M}$ is a matrix defined by
$$
\widetilde{m}_{ip} = \begin{cases}
m_{1p} + m_{2p}+ \cdots +m_{p-1,p} +m_{p,p} + m_{n,p} & \text{if $i=p$}, \\
m_{ip} & \text{if $p+1 \le i \le n-1$} \\
0 & \text{otherwise},
\end{cases}
$$
and $\widetilde{m}_{ij}=m_{ij}$ for $j \not=p$.
It is easy to see that $\widetilde{M}\in\mathbb{M}_{n,l}(p-1)$.
Then, using Lemma \ref{lem:reduction} below, we get
$$
E_{1p}^{m_{1p}} E_{2p}^{m_{2p}} \cdots E_{p-1,p}^{m_{p-1,p}} E_{np}^{m_{np}} \cdot 
\bm{e}^{\widetilde{M}} \equiv (\text{non-zero constant}) \times \bm{e}^M
\pmod{\bigoplus_{N \in \mathbb{M}_{n,l}(p-1)} \bC \cdot \bm{e}^{N}}.
$$
Therefore we have $M \in \widetilde{\mathbb{M}}_{n,l}$ by the induction hypothesis on $p$,
and hence $\mathbb{M}_{n,l}(p) \subset \tM nl$.
In particular, we get
$\mathbb{M}_{n,l}=\mathbb{M}_{n,l}(n-1)\subset \tM nl$,
which is the desired conclusion.
\end{proof}

In the proof, we have used the following lemma which is readily verified.
\begin{lem}\label{lem:reduction}
Suppose that $M\in\mathbb{M}_{n,l}(p-1)$ and $1\le i,k\le p-1$.
Then $M+R_k^{ip}=M-R_k^{pi}\in\mathbb{M}_{n,l}(p-1)$.
In particular,
\begin{align*}
E_{ip}^d\bm{e}^M \equiv m_{pp}(m_{pp}-1)\dotsb(m_{pp}-d+1)\bm{e}^{M-dR_p^{pi}}
\pmod{\bigoplus_{N \in \mathbb{M}_{n,l}(p-1)} \bC \cdot \bm{e}^{N}}
\end{align*}
holds for $d\ge1$.
\qed
\end{lem}

The linear map $\rho_{\gl_n}: \U \to \mr{End}_{\bC}(\mcal{A}(\mr{Mat}_n))$ defined by
$$
(\rho_{\gl_n} (E_{ij}) \cdot f) (X)= \sum_{k=1}^n x_{ik} \frac{\partial f}{\partial x_{jk}} (X) 
\qquad (1 \le i,j \le n, \quad f \in \mcal{A}(\mr{Mat}_n))
$$
determines a representation of $\U$ on $\mcal{A}(\mr{Mat}_n)$.
We abbreviate $\rho_{\gl_n}(E_{ij})$ as $E_{ij}$.

For each $M \in \mbb{M}_{n,l}$, put $X^M= \prod_{i,j=1}^n x_{ij}^{m_{ij}}$.
The action of $\U$ to $X^M$ is given by
$$
E_{pq} \cdot X^M = \sum_{k=1}^n m_{qk} X^{M+R_{k}^{pq}} \qquad (1 \le p,q \le n).
$$
Combining this with \eqref{eq:ActionSymTensor} and Lemma \ref{lem:Cyclic},
we see that the linear map 
$\bm{e}^M \mapsto X^M \ (M \in \mathbb{M}_{n,l})$ gives 
the isomorphism 
\begin{equation} \label{eq:Alpha=0} 
(\st^l (\bC^n))^{\otimes n} \cong \bigoplus_{M \in \mathbb{M}_{n,l}} \bC \cdot X^M 
=\U \cdot x_{11}^l x_{22}^l \cdots x_{nn}^l
\subset 
\mcal{A}(\mr{Mat}_n).
\end{equation}

\section{The cyclic modules $\U \cdot \dete{\alpha}(X)^l$}
%

\subsection{$\alpha$-determinants and intertwiners}

Let $\alpha$ be a complex number.
We consider the cyclic module $\U \cdot \dete{\alpha}(X)^l$
for a positive integer $l$.

When $\alpha=0$,
we have $\dete{0}(X) = x_{11} x_{22} \cdots x_{nn}$.
From \eqref{eq:Alpha=0} we obtain the irreducible decomposition
\begin{equation}
\U \cdot \dete{0}(X)^l \cong (\st^l(\bC^n))^{\otimes n} \cong
\bigoplus_{\lambda \vdash ln} (\mcal{M}_n^\lambda)^{\oplus K_{\lambda (l^n)}}.
\end{equation}
In general,
the module $\U \cdot \dete{\alpha}(X)^l$ is a submodule of 
$\U \cdot \dete{0}(X)^l$ because  
$\dete{\alpha}(X)^l \in \bigoplus_{M \in \mathbb{M}_{n,l}} \bC \cdot X^M
=\U \cdot \dete{0}(X)^l$.
Therefore we have

\begin{thm}\label{thm:decomposition}
It holds that
$$
\U \cdot \dete{\alpha}(X)^l \cong 
\bigoplus_{\lambda \vdash ln} (\mcal{M}_n^\lambda)^{\oplus m_{n,l}^{\lambda}(\alpha)},
$$
where $m_{n,l}^\lambda(\alpha)$ is a nonnegative integer at most $K_{\lambda (l^n)}$
and $m_{n,l}^\lambda(0)=K_{\lambda (l^n)}$. \qed
\end{thm}

In order to obtain further properties of the multiplicities $m_{n,l}^\lambda(\alpha)$,
we construct a $\U$-intertwiner from $(\st^{l}(\bC^n))^{\otimes n}$
to $\U \cdot \dete{\alpha}(X)^l$ explicitly for each $\alpha$.

For a sequence $(k_1,\dots, k_n ) \in [n]^{\times n}$, define
$$
D^{(\alpha)} (k_1,\dots, k_n) = \dete{\alpha} 
\begin{pmatrix} x_{k_1 1} & x_{k_1 2} & \ldots & x_{k_1 n} \\
x_{k_2 1} & x_{k_2 2} & \ldots & x_{k_2 n} \\
\vdots & \vdots & \ddots & \vdots \\
x_{k_n 1} & x_{k_n 2} & \ldots & x_{k_n n} \end{pmatrix}.
$$
For a matrix $N \in \mathbb{M}_{n,1}$,
there exists some $(k_1,\dots, k_n ) \in [n]^{\times n}$
such that
$N=(\delta_{i,k_j})_{1\le i,j\le n}$.
Then we let
$$
D^{(\alpha)} (N) = D^{(\alpha)}(k_1,\dots,k_n).
$$

Let $M=(m_{ij})_{1 \le i,j \le n} \in \mathbb{M}_{n,l}$.
A sequence $(M_1,\dots,M_l)\in(\mathbb{M}_{n,1})^{\times l}$
is called a \emph{partition} of $M$ and denoted by
$(M_1,\dots,M_l)\Vdash M$ if $M_1+\dots+M_l=M$.
We also put $M!=\prod_{i,j=1}^nm_{ij}!$.
For instance, $(lI_n)!=l!^n$.
Now we define 
the element $D^{(\alpha)}(M) \in \mcal{A}(\mr{Mat}_n)$ by
\begin{equation}
D^{(\alpha)}(M)=\frac{M!}{(lI_n)!}
\sum_{(M_1,\dots,M_l)\Vdash M}
D^{(\alpha)}(M_1) D^{(\alpha)}(M_2) \cdots D^{(\alpha)}(M_l),
\end{equation}
where the sum runs over all partitions of $M$.
It is clear that
$D^{(\alpha)}(lI_n) =\dete{\alpha}(X)^l$.

\begin{example}
\begin{align*}
D^{(\alpha)}\!\begin{pmatrix} 2 & 1 \\ 0 & 1 \end{pmatrix}
&= D^{(\alpha)}\!\begin{pmatrix} 1 & 1 \\ 0 & 0 \end{pmatrix}
D^{(\alpha)}\!\begin{pmatrix} 1 & 0 \\ 0 & 1 \end{pmatrix}
=D^{(\alpha)}(1,1) D^{(\alpha)}(1,2). \\
D^{(\alpha)}\!\begin{pmatrix} 2 & 1 \\ 1 & 2 \end{pmatrix}
&= \frac{(2!)^2}{(3!)^2}
\ckakko{6 D^{(\alpha)}\!\begin{pmatrix} 1 & 1 \\ 0 & 0 \end{pmatrix}
D^{(\alpha)}\!\begin{pmatrix} 1 & 0 \\ 0 & 1 \end{pmatrix}
D^{(\alpha)}\!\begin{pmatrix} 0 & 0 \\ 1 & 1 \end{pmatrix}
+ 3
D^{(\alpha)}\!\begin{pmatrix} 1 & 0 \\ 0 & 1 \end{pmatrix}^{\!2}
D^{(\alpha)}\!\begin{pmatrix} 0 & 1 \\ 1 & 0 \end{pmatrix}} \\
&= \frac{1}{3}
\(2 D^{(\alpha)}(1,1) D^{(\alpha)}(1,2) D^{(\alpha)}(2,2)+
 D^{(\alpha)}(1,2)^2 D^{(\alpha)}(2,1) \).
\end{align*}
\qed
\end{example}

Take $M=(m_{ij})_{1\le i,j\le n}\in\mathbb{M}_{n,l}$
and suppose that $m_{qk}>0$.
Then $M+R_k^{pq}\in\mathbb{M}_{n,l}$.
Let $(M_1,\dots,M_l)\Vdash M$ and
$(M_1',\dots,M_l')\Vdash M+R_k^{pq}$.
We write
$(M_1,\dots,M_l)\xrightarrow{p,q;\,k}(M_1',\dots,M_l')$
if there exists some $j$ such that
$M_i'=M_i+\delta_{ij}R_k^{pq}$.
We notice that
\begin{equation}\label{eq:number_of_adjacent_partitions}
\#\Bigl\{
(M_1',\dots,M_l')\Vdash M+R_k^{pq}\ \Big|\
(M_1,\dots,M_l)\xrightarrow{p,q;\,k}(M_1',\dots,M_l')
\Bigr\}=m_{qk}
\end{equation}
because $M_j+R_k^{pq}\in\mathbb{M}_{n,1}$
if and only if $(M_j)_{qk}=1$ so that
the number of such choices of $j$ is just
$m_{qk}=\sum_{j=1}^l(M_j)_{qk}$.
We also notice that
\begin{align*}
(M_1,\dots,M_l)\xrightarrow{p,q;\,k}(M_1',\dots,M_l')
\iff
(M_1',\dots,M_l')\xrightarrow{q,p;\,k}(M_1,\dots,M_l).
\end{align*}


The following fact is crucial.
\begin{prop} \label{prop:ActionAlphaDet}
For any $p, q \in [n]$ and  $M \in \mathbb{M}_{n,l}$, we have
\begin{equation} \label{eq:ActionAlphaDet}
E_{pq} \cdot D^{(\alpha)}(M) = \sum_{k=1}^n m_{qk} D^{(\alpha)} (M +R^{pq}_k). 
\end{equation}
\end{prop}

\begin{example}
\begin{align*}
E_{11} \cdot D^{(\alpha)}\!\begin{pmatrix} 2 & 1 \\ 1 & 2 \end{pmatrix} 
&=3 D^{(\alpha)}\!\begin{pmatrix} 2 & 1 \\ 1 & 2 \end{pmatrix}, 
& \quad
E_{12} \cdot D^{(\alpha)}\!\begin{pmatrix} 2 & 1 \\ 1 & 2 \end{pmatrix}
&=D^{(\alpha)}\!\begin{pmatrix} 3 & 1 \\ 0 & 2 \end{pmatrix}+%
2 D^{(\alpha)}\!\begin{pmatrix} 2 & 2 \\ 1 & 1 \end{pmatrix},  \\
E_{21} \cdot D^{(\alpha)}\!\begin{pmatrix} 2 & 1 \\ 1 & 2 \end{pmatrix}
&=2 D^{(\alpha)}\!\begin{pmatrix} 1 & 1 \\ 2 & 2 \end{pmatrix}+%
D^{(\alpha)}\!\begin{pmatrix} 2 & 0 \\ 1 & 3 \end{pmatrix},
& \quad
 E_{22} \cdot D^{(\alpha)}\!\begin{pmatrix} 2 & 1 \\ 1 & 2 \end{pmatrix} 
&=3 D^{(\alpha)}\!\begin{pmatrix} 2 & 1 \\ 1 & 2 \end{pmatrix}. 
\end{align*}
\qed
\end{example}

\begin{proof}[Proof of Proposition \ref{prop:ActionAlphaDet}]
First we notice that
we can verify the case where $l=1$ easily
(see Lemma 2.1 in \cite{MW}).
By using this result,
for any $M \in \mathbb{M}_{n,l}$,
we have
\begin{align*}
E_{pq} \cdot D^{(\alpha)}(M)
&=
\frac{M!}{(lI_n)!} 
\sum_{(M_1,\dots, M_l)\Vdash M} \sum_{j=1}^l
D^{(\alpha)}(M_1) \cdots (E_{pq} \cdot D^{(\alpha)}(M_j)) \cdots D^{(\alpha)}(M_l)  \\
&=\frac{M!}{(lI_n)!}
\sum_{k=1}^n \sum_{(M_1,\dots, M_l)\Vdash M}
\sum_{j=1}^l(M_j)_{qk} 
D^{(\alpha)}(M_1) \cdots  D^{(\alpha)}(M_j+R_{pq}^k) \cdots D^{(\alpha)}(M_l) \\
&=\frac{M!}{(lI_n)!}
\sum_{k=1}^n \sum_{(M_1,\dots, M_l)\Vdash M}
\sum_{\substack{(M_1',\dots,M_l')\Vdash M+R_k^{pq}\\
(M_1,\dots,M_l)\xrightarrow{p,q;\,k}(M_1',\dots,M_l')}}
D^{(\alpha)}(M_1') \cdots   D^{(\alpha)}(M_l')\\
&=\frac{M!}{(lI_n)!}
\sum_{k=1}^n \sum_{(M_1',\dots,M_l')\Vdash M+R_k^{pq}}
\sum_{\substack{(M_1,\dots, M_l)\Vdash M\\
(M_1',\dots,M_l')\xrightarrow{q,p;\,k}(M_1,\dots,M_l)}}
D^{(\alpha)}(M_1')\dots D^{(\alpha)}(M_l').
\end{align*}
By \eqref{eq:number_of_adjacent_partitions},
we see that
\begin{align*}
\sum_{\substack{(M_1,\dots, M_l)\Vdash M\\
(M_1',\dots,M_l')\xrightarrow{q,p;\,k}(M_1,\dots,M_l)}}1
&=\#\Bigl\{
(M_1,\dots,M_l)\Vdash (M+R_k^{pq})+R_k^{qp}\ \Big|\
(M_1',\dots,M_l')\xrightarrow{q,p;\,k}(M_1,\dots,M_l)
\Bigr\}\\
&=(M+R_k^{pq})_{pk}=m_{pk}+1.
\end{align*}
Hence it follows that
\begin{align*}
E_{pq} \cdot D^{(\alpha)}(M)
=& 
\sum_{k=1}^n (m_{pk}+1)\frac{M!}{(lI_n)!}
\sum_{(M'_1,\dots, M'_l)\Vdash M+R_k^{pq}}
D^{(\alpha)}(M_1') \cdots   D^{(\alpha)}(M_l')  \\
=& \sum_{k=1}^n m_{qk} D^{(\alpha)}(M+R^{pq}_k)
\end{align*}
since $(m_{pk}+1)M!=m_{qk}(M+R_k^{pq})!$ if $m_{qk}>0$.
Thus we have proved \eqref{eq:ActionAlphaDet}.
\end{proof}

Now we give an explicit intertwiner
from $(\st^l(\bC^n))^{\otimes n}$
to $\U \cdot \dete{\alpha}(X)^l$.
The following proposition is a generalization 
of Lemma 2.3 and Proposition 2.4 in \cite{MW} for the case where $l=1$.

\begin{prop} \label{prop:intertwiner}
We have
$$
\U \cdot \dete{\alpha}(X)^l =
\sum_{M \in \mathbb{M}_{n,l}} \bC \cdot D^{(\alpha)}(M).
$$
Furthermore, the linear map $\Phi^{(\alpha)}$ 
determined by
$$
\Phi^{(\alpha)}(\bm{e}^M) = D^{(\alpha)}(M), \qquad M \in \mathbb{M}_{n,l},
$$
gives a surjective $\U$-intertwiner 
from $(\st^l(\bC^n))^{\otimes n}$
to $\U \cdot \dete{\alpha}(X)^l$.
\end{prop}

\begin{proof}
From Proposition \ref{prop:ActionAlphaDet}, the space
$\sum_{M \in \mathbb{M}_{n,l}} \bC \cdot D^{(\alpha)}(M)$ is invariant under 
the action of $\U$.
Since $D^{(\alpha)}(l I_n )= \dete{\alpha}(X)^l$,
the space $\U \cdot \dete{\alpha}(X)^l$ is a submodule of 
$\sum_{M \in \mathbb{M}_{n,l}} \bC \cdot D^{(\alpha)}(M)$.
Furthermore, by \eqref{eq:ActionSymTensor} and Proposition \ref{prop:ActionAlphaDet},
the linear map $\Phi^{(\alpha)}$ 
determined by
$$
\Phi^{(\alpha)}(\bm{e}^M) = D^{(\alpha)}(M), \qquad M \in \mathbb{M}_{n,l},
$$
gives a surjective $\U$-intertwiner 
from $(\st^l(\bC^n))^{\otimes n}$
to $\sum_{M \in \mathbb{M}_{n,l}} \bC \cdot D^{(\alpha)}(M)$.
It follows from Lemma \ref{lem:Cyclic} that 
$$
\sum_{M \in \mathbb{M}_{n,l}} \bC \cdot D^{(\alpha)}(M)
= \sum_{M \in \mathbb{M}_{n,l}} \bC \cdot \Phi^{(\alpha)}(\bm{e}^M)
 \subset \U \cdot \Phi^{(\alpha)}
(\bm{e}^{l I_n}) = \U \cdot \dete{\alpha}(X)^l
$$ 
as we desired.
\end{proof}

\subsection{Transition matrices}

We show that the multiplicity $m_{n,l}^{\lambda}(\alpha)$
in Theorem \ref{thm:decomposition}
is described as a rank of a certain matrix
for each highest weight $\lambda$.

The module $(\st^l(\bC^n))^{\otimes n}$ is decomposed in the form
$$
(\st^l(\bC^n))^{\otimes n}  = \bigoplus_{\lambda \vdash nl} 
\bigoplus_{i=1}^{K_{\lambda (l^n)}} \U \cdot v_{i}^\lambda.
$$
Here $v_i^\lambda$ ($i=1,\dots,K_{\lambda(l^n)}$)
are highest weight vectors corresponding to the weight $\lambda$.
Under the isomorphism $\Phi^{(0)}$ and 
surjective intertwiner $\Phi^{(\alpha)}$, we see that
\begin{align*}
\U \cdot \dete{0}(X)^l =&
\bigoplus_{\lambda \vdash nl} 
\bigoplus_{i=1}^{K_{\lambda (l^n)}} \U \cdot \Phi^{(0)} (v_{i}^\lambda), \\
\U \cdot \dete{\alpha}(X)^l =&
\bigoplus_{\lambda \vdash nl} 
\sum_{i=1}^{K_{\lambda (l^n)}} \U \cdot \Phi^{(\alpha)} (v_{i}^\lambda).
\end{align*}
Since  $\Phi^{(\alpha)}(v_i^\lambda)$ is the highest weight vector
unless it vanishes, there exists a matrix 
$\tramat\lambda(\alpha)=((\tramat\lambda(\alpha))_{ij}) $ of size $K_{\lambda (l^n)}$
such that
\begin{equation} \label{eq:TransitionMatrix}
\Phi^{(\alpha)}(v_j^\lambda) = \sum_{i=1}^{K_{\lambda (l^n)}} (\tramat\lambda(\alpha))_{ij}
\Phi^{(0)}(v_i^\lambda)
\end{equation}
for each $j$.
We call the matrix $\tramat\lambda(\alpha)$ the {\it transition matrix}.
We notice that the definition of $\tramat\lambda(\alpha)$ is
\emph{dependent} on the choice of vectors $v_1^\lambda, \dots, v_{K_{\lambda
(l^n)}}^\lambda$
but $\tramat\lambda(\alpha)$ is \emph{uniquely determined up to conjugacy}.
By definition, its entries belong to $\bQ[\alpha]$.
We now obtain the

\begin{thm} \label{thm:matrixF}
For each $\alpha \in \bC$ and $\lambda \vdash nl$,
the multiplicity $m_{n,l}^\lambda(\alpha)$
in Theorem \ref{thm:decomposition} is
equal to the rank of the matrix
$\tramat\lambda(\alpha)$ defined via \eqref{eq:TransitionMatrix}.
Namely, the irreducible decomposition of the cyclic module
$\U\cdot\adet(X)^l$ is given by
\begin{align}
\U\cdot\adet(X)^l\cong\bigoplus_{\substack{\lambda\vdash nl\\\len\lambda\le n}}
(\GLmod \lambda)^{\oplus\rank\tramat\lambda(\alpha)}.
\end{align}
\qed
\end{thm}

We need to obtain an explicit expression
of the matrix $\tramat\lambda(\alpha)$
to evaluate the multiplicity $m_{n,l}^\lambda(\alpha)$.
When $n=2$, we show that
the matrix $\tramat[2,l]\lambda(\alpha)$
is of size $1$ and given explicitly by a hypergeometric polynomial.
See the next section for the detailed discussion for this case.
In general,
it is not easy to calculate
the matrix $\tramat\lambda(\alpha)$,
and we have no effective method to evaluate the multiplicity
$m_{n,l}^\lambda(\alpha)$.

We give several examples of
an explicit calculation of transition matrices
for the highest weights with special types.

\begin{example}
If $l=1$, then we have
$\tramat[n,1]\lambda(\alpha)=f_\lambda(\alpha)I$
for any partition $\lambda$ of $n$,
where $f_\lambda(\alpha)$ is defined in 
\eqref{Polynomial_f}.
See Corollary 3.4 in \cite{MW}. \qed
\end{example}

\begin{example}
For $\lambda=(nl)$, the Kostka number $K_{\lambda (l^n)}$ is equal to $1$. 
The vector $v^{(nl)} =\bm{e}_1^l \otimes \bm{e}_{1}^l \otimes \cdots \otimes \bm{e}_1^l$
is the highest weight vector with the highest weight $(nl)$.
By Proposition \ref{prop:intertwiner} we have
$$
\Phi^{(\alpha)}(v^{(nl)})
=D^{(\alpha)}(1,1,\dots,1)^l
=\Bigl\{\prod_{j=1}^{n-1} (1+j\alpha)  x_{11} x_{22} \cdots x_{nn}\Bigr\}^{\!l}
= \prod_{j=1}^{n-1} (1+j\alpha)^l \cdot \Phi^{(0)}(v^{(nl)})
$$
and hence
\begin{equation*}
\tramat{(nl)}(\alpha)=\prod_{j=1}^{n-1} (1+j \alpha)^l. 
\end{equation*}
\qed
\end{example}

\begin{example}\label{ex:standard_case}
For $\lambda=(nl-1,1)$, the Kostka number $K_{\lambda (l^n)}$ is equal to $n-1$. 
Put 
$$
w_i= \bm{e}_1^l \otimes \bm{e}_1^l \otimes \cdots \otimes 
\stackrel{\text{$i$-th}}{\bm{e}_1^{l-1} \bm{e}_2} \otimes \cdots 
\otimes \bm{e}_1^l
$$
for each $1 \le i \le n$.
Then $v_i^{(nl-1,1)}= w_i-w_{i+1}$ $(1 \le i \le n-1)$ are
linearly independent highest weight vectors
corresponding to the weight $(nl-1,1)$.
It is easy to see that
$$
\Phi^{(\alpha)}(v_i^{(nl-1,1)}) = (1-\alpha)(1+(n-1)\alpha)^{l-1}  
\prod_{j=1}^{n-2} (1+j\alpha)^l \cdot \Phi^{(0)} (v_i^{(nl-1,1)})
$$
which readily implies
$$
\tramat{(nl-1,1)}(\alpha)=(1-\alpha)(1+(n-1)\alpha)^{l-1}  
\prod_{j=1}^{n-2} (1+j\alpha)^l \cdot I_{n-1}. 
$$
\qed
\end{example}

%
\section{$\U[2]$-cyclic modules and Jacobi polynomials}
%

In this section, we study the case where $n=2$.
The transition matrix $\tramat[2,l]\lambda(\alpha)$ is of size $1$
and explicitly given by a hypergeometric polynomial in $\alpha$
which is in fact the Jacobi polynomial.
Moreover, we see that these Jacobi polynomials are unitary.

\subsection{Explicit irreducible decomposition of $\U[2] \cdot \dete{\alpha}(X)^l$}

For a non-negative integer $n$,
complex numbers $b$ and $c$ such that
$c \not= -1,-2,\dots,-n+1$,
let $F(-n,b,c;x)$ be the Gaussian hypergeometric polynomial
$$
F(-n,b,c;x)= 1+\sum_{k=1}^n \frac{(-n)_k (b)_k}{(c)_k} \frac{x^k}{k!}.
$$
Here $(a)_k$ stands for the Pochhammer symbol $(a)_k= a(a+1) \cdots (a+k-1)$.
For any partition $\lambda$ of $2l$ with length $\le 2$,
we have $K_{\lambda (l^2)}=1$,
whence $\tramat[2,l]\lambda(\alpha)$ is a scalar.

\begin{thm}\label{thm:CaseTwo}
For non-negative integers $l$ and $s$ such that $0 \le s \le l$, we have
\begin{equation}
\tramat[2,l]{(2l-s,s)}(\alpha) = (1+\alpha)^{l-s} G_s^l(\alpha),
\end{equation}
where $G_{n}^\gamma (x)$ is the polynomial given by
$$
G_{n}^\gamma (x) = F(-n,\gamma-n+1,-\gamma;-x).
$$
\end{thm}
By the hypergeometric 
differential equation satisfied by $G_{n}^\gamma (x)$, 
the explicit form
of $\tramat[2,l]{(2l-s,s)}(\alpha)$ given in Theorem \ref{thm:CaseTwo} shows that
 $\tramat[2,l]{(2l-s,s)}(\alpha)$ satisfies 
the following singly confluent Heun differential equation 
(see \cite{SL}).
\begin{cor}\label{cor:Heun}
The polynomial $f(x)=\tramat[2,l]{(2l-s,s)}(-x)$
satisfies the differential equation
\begin{equation}\label{eq:deq_for_F}
\ckakko{\frac{d^2}{dx^2}+\kakko{\frac2{x-1}-\frac{l}x}\frac{d}{dx}%
+\frac{s-(l-s)^2-x}{x(x-1)^2}}f(x)=0.
\end{equation}
\end{cor}

\begin{remark}\label{remark:Heun}
Since it seems difficult at present
to obtain the transition matrices explicitly in general,
we are naturally lead to the following questions:
Can one obtain the equation \eqref{eq:deq_for_F} directly
by investigating (a certain structure of)
the cyclic module $\U[2]\cdot\adet(X)^l$ itself?
If it is possible, is the derivation of the differential equation
generalized to the cases where $n\ge3$?
\end{remark}

The roots of the polynomial $G^\gamma_n(x)$ satisfy the following property.

\begin{prop} \label{prop:unitarity}
For a real number $\gamma$ such that $\gamma \ge n$, the polynomial
$G^\gamma_n(x)$ is unitary, i.e., 
every root of $G^\gamma_n(x)$ is
on the unit circle $\bT=\{z\in\bC \ | \ |z|=1\}$.
Furthermore, 
$G_{2n}^\gamma(1)\ne0$ and $G^\gamma_{2n+1}(1)=0$
for any nonnegative integer $n$.
\end{prop}

Therefore we obtain the following irreducible decomposition 
from Theorem \ref{thm:decomposition} and Proposition \ref{thm:matrixF}.

\begin{cor}
For any $\alpha \in \bC \setminus \bT$, we have
$$
\U[2] \cdot \dete{\alpha}(X)^l \cong (\st^l(\bC^2))^{\otimes 2}
\cong \bigoplus_{s=0}^l \mcal{M}_2^{(2l-s,s)}.
$$
For $\alpha=\pm 1$, we have
\begin{align*}
\U[2]\cdot\per(X)^l
&\cong \bigoplus_{j=0}^{\floor{l/2}} \mcal{M}_2^{(2l-2j,2j)}
\cong \st^l (\st^2 (\bC^2)),\\
\U[2] \cdot \det(X)^l
&=\bC \cdot \det(X)^l
\cong \mcal{M}_2^{(l,l)}.
\end{align*}
\end{cor}

\subsection{Proof of Theorem \ref{thm:CaseTwo}}

The highest weight vector associated with the highest weight $(2l-s,s)$
in the module $(\st^l(\bC^2))^{\otimes 2}$ ($\cong
\bigoplus_{s=0}^l \mcal{M}_2^{(2l-s,s)}$)
is given by
$$
v^{(2l-s,s)} = \sum_{j=0}^s (-1)^j \binom{s}{j} \bm{e}_1^{l-j} \bm{e}_2^j \otimes 
\bm{e}_1^{l-s+j} \bm{e}_2^{s-j}.
$$
The image of this under $\Phi^{(\alpha)}$ is 
\begin{equation} \label{eq:PhiV1}
\Phi^{(\alpha)}(v^{(2l-s,s)}) = 
\sum_{j=0}^s (-1)^j \binom{s}{j} D^{(\alpha)}\!
\begin{pmatrix} l-j & l-s+j \\ j & s-j \end{pmatrix}.
\end{equation}

\begin{lem} \label{lem:DMpq}
For $0 \le p \le q \le l$,  we have
\begin{align*}
&D^{(\alpha)}\!\begin{pmatrix} l-p & l-q \\ p & q \end{pmatrix} \\
=&\binom{l}{q}^{-1}
\sum_{r=0}^{\min\{p,l-q\}} \binom{l-p}{q-p+r} \binom{p}{r}
 D^{(\alpha)}(1,1)^{l-q-r}
D^{(\alpha)}(1,2)^{q-p+r}D^{(\alpha)}(2,1)^{r} D^{(\alpha)}(2,2)^{p-r}.
\end{align*}
\end{lem}

\begin{proof}
Sequences $(M_1,\dots,M_l) \in (\mathbb{M}_{2,1})^{\times l}$ satisfying
$M_1+ \cdots+M_l=\begin{smatrix} l-p & l-q \\ p & q \end{smatrix}$
are permutations of 
$$
\overbrace{
\begin{pmatrix} 1 & 1 \\ 0 & 0 \end{pmatrix}, \dots, 
\begin{pmatrix} 1 & 1 \\ 0 & 0 \end{pmatrix}}^{l-q-r}, \quad
\overbrace{
\begin{pmatrix} 1 & 0 \\ 0 & 1 \end{pmatrix}, \dots, 
\begin{pmatrix} 1 & 0 \\ 0 & 1 \end{pmatrix}}^{q-p+r}, \quad
\overbrace{
\begin{pmatrix} 0 & 1 \\ 1 & 0 \end{pmatrix}, \dots, 
\begin{pmatrix} 0 & 1 \\ 1 & 0 \end{pmatrix}}^{r}, \quad
\overbrace{
\begin{pmatrix} 0 & 0 \\ 1 & 1 \end{pmatrix}, \dots, 
\begin{pmatrix} 0 & 0 \\ 1 & 1 \end{pmatrix}}^{p-r},
$$
where $r$ runs over  $0,1,\dots, \min \{p,l-q\}$.
Since the number of such sequences is $(l!)/\{(l-q-r)! \, (q-p+r)! \, (p-r)!\}$,
we have
\begin{align*}
D^{(\alpha)}\begin{pmatrix} l-p & l-q \\ p & q \end{pmatrix}
=& \frac{(l-p)! \, (l-q)! \, p! \, q!}{(l!)^2} \sum_{r=0}^{\min\{p,l-q\}}
\frac{l!}{(l-q-r)! \, (q-p+r)! \, r! \, (p-r)!}  \\
& \times D^{(\alpha)}(1,1)^{l-q-r}
D^{(\alpha)}(1,2)^{q-p+r} D^{(\alpha)}(2,1)^{r} D^{(\alpha)}(2,2)^{p-r}.
\end{align*}
This completes the proof.
\end{proof}

The polynomial $\tramat[2,l]{(2l-s,s)}(\alpha)$ is determined by the identity
$\Phi^{(\alpha)}(v^{(2l-s,s)})=\tramat[2,l]{(2l-s,s)}(\alpha)\Phi^{(0)}(v^{(2l-s,s)})$.
By comparing the coefficients of $x_{11}^l x_{12}^{l-s}x_{22}^s$ in the both sides, 
we see that
$$
\tramat[2,l]{(2l-s,s)}(\alpha) = [x_{11}^l x_{12}^{l-s} x_{22}^s] \,
\Phi^{(\alpha)}(v^{(2l-s,s)}).
$$
Here $[x_{11}^ax_{12}^bx_{21}^cx_{22}^d]f(x_{11},x_{12},x_{21},x_{22})$
stands for the coefficient of $x_{11}^ax_{12}^bx_{21}^cx_{22}^d$ in $f(x_{11},x_{12},x_{21},x_{22})$.
By using Lemma \ref{lem:DMpq} together with
\begin{align*}
&& D^{(\alpha)}(1,1)&=(1+\alpha) x_{11} x_{12}, &
D^{(\alpha)}(1,2)&=x_{11}x_{22}+ \alpha x_{21}x_{12}, &\\
&& D^{(\alpha)}(2,1)&=\alpha x_{11}x_{22}+  x_{21}x_{12}, & 
D^{(\alpha)}(2,2)&=(1+\alpha) x_{21}x_{22},  &
\end{align*}
we have
\begin{align*}
&\eqspace[x_{11}^l x_{12}^{l-s} x_{22}^s] 
D^{(\alpha)}\begin{pmatrix} l-j & l-s+j \\ j & s-j \end{pmatrix} \\
&=\binom{l}{s-j}^{-1}
\binom{l-j}{s-j}
\kakko{[(x_{11} x_{12})^{l-s}] D^{(\alpha)}(1,1)^{l-s}}\cdot 
\kakko{[(x_{11}x_{22})^{s}] D^{(\alpha)}(1,2)^{s-j}D^{(\alpha)}(2,1)^{j}}\\
&=\frac{(l-j)! \, (l-s+j)!}{l! \, (l-s)!}(1+\alpha)^{l-s} \alpha^j 
\end{align*}
for $0\le j<s/2$.
We can check that this identity holds for any $0 \le j \le s$ in a similar way.
Hence it follows from \eqref{eq:PhiV1} that 
\begin{align*}
\tramat[2,l]{(2l-s,s)}(\alpha)
&=\frac{s!}{l! \, (l-s)!} (1+\alpha)^{l-s}
\sum_{j=0}^s \frac{(l-j)! \, (l-s+j)!}{(s-j)!} \cdot \frac{(-\alpha)^j}{j!}\\
&=(1+\alpha)^{l-s}
\sum_{j=0}^s \frac{s!\,(l-j)! \, (l-s+j)!}{l! \, (l-s)!\,(s-j)!} \cdot \frac{(-\alpha)^j}{j!}\\
&=(1+\alpha)^{l-s}F(-s,l-s+1,-l;-\alpha).
\end{align*}
Thus we have proved Theorem \ref{thm:CaseTwo}.

\subsection{Proof of Proposition \ref{prop:unitarity}}

Let $\gamma$ be a positive real number
and $n$ a non-negative integer such that $\gamma \ge n$.
We prove the unitarity of the polynomial $G_n^\gamma(x)$
by the property of the Jacobi polynomial (see e.g. \cite{Szego})
$$
P_n^{(\alpha,\beta)}(x)=
\binom{n+\alpha}{n}F\!\kakko{-n,n+\alpha+\beta+1, \alpha+1; \frac{1-x}{2}}.
$$
We see by definition that 
\begin{equation} \label{eq:JacobiG}
G_n^\gamma(x) = \binom{n-\gamma-1}{n}^{-1} P_n^{(-\gamma-1,2\gamma-2n+1)}(1+2x).
\end{equation}
We recall the following formulas
((4.1.3), (4.22.1), and (4.1.5) in \cite{Szego}).
\begin{lem} 
For any $\alpha, \beta \in \bC$ and non-negative integer $n$, the following formulas hold.
\begin{align}
P_{n}^{(\alpha,\beta)}(x)=&(-1)^n P_n^{(\beta,\alpha)}(-x), \label{eq:JacobiFE} \\
P_n^{(\alpha,\beta)}(x)=& \(\frac{1-x}{2}\)^n P_n^{(-2n-\alpha-\beta-1,\beta)} \(\frac{x+3}{x-1}\), 
\label{eq:JacobiFE2} \\
P_{2n}^{(\alpha,\beta)}(x)=& \frac{(-1)^n n!}{(2n)!} (\alpha+n+1)_n P_n^{(-\frac{1}{2},\alpha)}
(1-2x^2), \label{eq:JacobiFEeven} \\
P_{2n+1}^{(\alpha,\beta)}(x)=& \frac{(-1)^n n!}{(2n+1)!} 
(\alpha+n+1)_{n+1} x P_n^{(\frac{1}{2},\alpha)}
(1-2x^2). \label{eq:JacobiFEodd}
\end{align}
\qed
\end{lem}

From \eqref{eq:JacobiG} and \eqref{eq:JacobiFE},
we have
$$
G^\gamma_n(x)= \binom{n-\gamma-1}{n}^{-1} (-1)^n P_n^{(2\gamma-2n+1,-\gamma-1)}(-1-2x).
$$
By \eqref{eq:JacobiFE2}, it follows
$$
G^\gamma_n(x)= 
\binom{n-\gamma-1}{n}^{-1} (-1)^n (1+x)^n P_n^{(-\gamma-1,-\gamma-1)} \( \frac{x-1}{x+1} \).
$$
Applying \eqref{eq:JacobiFEeven} and \eqref{eq:JacobiFEodd} to this expression,
we obtain the following lemma.

\begin{lem} \label{Lem:GJacobiPexpression}
It holds that
\begin{align*}
G^\gamma_{2m}(x)
=& \frac{(-1)^m m!}{(-\gamma)_m} (x+1)^{2m} P_m^{(-1/2,-\gamma-1)}\(1-2\(\frac{x-1}{x+1}\)^2\), \\
G^\gamma_{2m+1}(x)
=& \frac{(-1)^{m+1} m!}{(-\gamma)_m} (x-1)(x+1)^{2m} 
P_m^{(1/2,-\gamma-1)}\(1-2\(\frac{x-1}{x+1}\)^2\).
\end{align*}
In particular, $G^\gamma_{2m}(1)\not= 0$ and $G^\gamma_{2m+1}(1)= 0$.
\qed
\end{lem}

In general,
the distribution of the roots of Jacobi polynomials
are described as follows.

\begin{lem}[Theorem 6.72 in \cite{Szego}]
For any real number $u$, let $E(u)$ be the Klein symbol, i.e.,
$$
E(u) = \begin{cases} 
u-1, & \text{if $u$ is an integer and $u \ge 0$}, \\
\floor u, & \text{if $u$ is not an integer and $u \ge 0$}, \\
0 & \text{if $u < 0$}.
\end{cases}
$$
Let $\alpha$ and $\beta$ be complex numbers and $n$ a non-negative integer.
Assume $\prod_{k=1}^n (\alpha+k)(\beta+k)(n+\alpha+\beta+k) \not=0$.
Define three numbers $X$, $Y$, and $Z$ by 
\begin{align*}
X=& E\( \frac{1}{2}( |2n+\alpha+\beta+1| - |\alpha| - |\beta| +1)\), \\
Y=& E\( \frac{1}{2}( -|2n+\alpha+\beta+1| + |\alpha| - |\beta| +1)\), \\
Z=& E\( \frac{1}{2}( -|2n+\alpha+\beta+1| - |\alpha| + |\beta| +1)\).
\end{align*}
Then, if we denote by $N(I)$ the number of roots of $P_n^{(\alpha,\beta)}(x)$
on an interval $I \subset \bR$, we have
\begin{align*}
N((-1,1))=& \begin{cases} 2 \floor{(X+1)/2}, 
& \text{if $(-1)^n \binom{n+\alpha}{n} \binom{n+\beta}{n} >0$}, \\
2 \floor{X/2}+1, 
& \text{if $(-1)^n \binom{n+\alpha}{n} \binom{n+\beta}{n} <0$}, 
\end{cases} \\
N((-\infty,-1))=& \begin{cases} 2 \floor{(Y+1)/2},
& \text{if $\binom{2n+\alpha+\beta}{n} \binom{n+\beta}{n} >0$}, \\
2 \floor{Y/2}+1, 
& \text{if $\binom{2n+\alpha+\beta}{n} \binom{n+\beta}{n} <0$}, 
\end{cases} \\
N((1,\infty))=& \begin{cases} 2 \floor{(Z+1)/2}, 
& \text{if $\binom{2n+\alpha+\beta}{n} \binom{n+\alpha}{n} >0$}, \\
2 \floor{Z/2}+1, 
& \text{if $\binom{2n+\alpha+\beta}{n} \binom{n+\alpha}{n} <0$}.
\end{cases}
\end{align*}
\qed
\end{lem}

Let us set $\alpha=-1/2$, $\beta=-\gamma-1$ and $n=2m$
(resp. $\alpha=1/2$, $\beta=-\gamma-1$ and $n=2m+1$)
in the lemma above and assume that $\gamma \ge n$.
It follows that $X=Y=0$ and $Z=m$,
from which we have
$N((-1,1))=N((-\infty,-1))=0$ and $N((1,\infty))=m$.
Since the degree of the polynomial
$P_m^{(-1/2,-\gamma-1)}(x)$ (resp. $P_m^{(1/2,-\gamma-1)}(x)$ ) is $m$,
all roots of $P_m^{(-1/2,-\gamma-1)}(x)$ (resp. $P_m^{(1/2,-\gamma-1)}(x)$)
belong to the interval $(1, \infty)$.
Therefore it follows from Lemma \ref{Lem:GJacobiPexpression} that
$$
a \in \bC,\quad G_n^\gamma(a) =0
\quad \Longrightarrow \quad
\frac{a-1}{a+1}\in i\bR
\quad \Longrightarrow \quad
|a|=1.
$$
This completes the proof of Proposition \ref{prop:unitarity}.

%
\section{Several remarks on the future study}
%

We give here several comments for the future study.

\subsection{Permanent cases}

When $\alpha=-1$,
$\adet[-1](X)$ is just the ordinary determinant
and we can easily see that
the cyclic module $\U \cdot \dete{-1}(X)^l$
is isomorphic to $\mcal{M}_n^{(l^n)}$.
However, in the case where $\alpha=1$,
we have not obtained the irreducible decomposition of
the cyclic module $\U\cdot\adet[1](X)^l$
generated by the permanent
$\per(X)=\adet[1](X)$.
Actually, only we can give here is the following conjecture.

\begin{conj} \label{Conj:Permanent}
$\U \cdot \dete{1}(X)^l \cong \st^l (\st^n(\bC^n))$.
\end{conj}

This claim is equivalent to the assertion that
the character of $\U \cdot \dete{1}(X)^l$ is
given by the plethysm $h_l \circ h_n$.
(For the definition of the plethysm for symmetric functions,
see \cite[\S I-8]{Mac}).
We have already verified
this conjecture in the following cases:
(i) $l=1$ (see \cite{MW}),
(ii) $n=1,2$ (see the previous section),
(iii) $n=3$ and $l=2$
(see Example \ref{ex:perm_for_(3,2)} below).

\begin{example}\label{ex:perm_for_(3,2)}
Let $n=3$ and $l=2$.
If we take a suitable highest weight vectors and 
employ a similar calculation in the proof of Theorem \ref{thm:CaseTwo}, 
we have
\begin{align*}
\tramat[3,2]{(6)}(\alpha)&=(1+\alpha)^2 (1+2\alpha)^2, \\
\tramat[3,2]{(5,1)}(\alpha)&=(1-\alpha)(1+\alpha)^2 (1+2\alpha) I_{2}, \\
\tramat[3,2]{(4,2)}(\alpha)&=(1+\alpha)^2 \cdot 
\diag \( 2(1-\alpha), 2(1-\alpha),  2-2\alpha+3\alpha^2 \), \\
\tramat[3,2]{(4,1,1)}(\alpha)&=\frac{1}{2}(1-\alpha)(1+\alpha)(2-5\alpha^2), \\
\tramat[3,2]{(3,3)}(\alpha)&=(1-\alpha)^2 (1+\alpha^2), \\
\tramat[3,2]{(3,2,1)} (\alpha)&=\frac{1}{4} (1-\alpha)(1+\alpha)(4-6\alpha+5\alpha^2) I_{2},\\
\tramat[3,2]{(2,2,2)}(\alpha)&=\frac{1}{2} (1-\alpha)^2 (2-2\alpha+5\alpha^2).
\end{align*}
In particular, when $\alpha=1$ we see that
\begin{align*}
m_{3,2}^\lambda(1)=\begin{cases}
1 & \lambda=(6),\,(4,2),\\
0 & \text{otherwise}
\end{cases}
\end{align*}
and hence it follows that
\begin{align*}
\U[3]\cdot\per(X)^2\cong\GLmod[3]{(6)}\oplus\GLmod[3]{(4,2)}
\cong\st^2(\st^3(\C^3))
\end{align*}
which agrees with our conjecture.
\end{example}

\subsection{Complex powers of $\alpha$-determinants} \label{Subsection:ComplexCase}

An appropriate reformulation of a setting is necessary
to study
``the cyclic module $\U \cdot \dete{\alpha}(X)^s$''
with a \emph{complex} number $s$.
Here we introduce a suitable space in which
we can treat such cyclic modules.

We take a $\U$-submodule
$$
\mathrm{ML}_n^{\bullet} = 
\left\{F_1 \cdots F_k \ | \ k \ge 0, \ 
F_i \in \bigoplus_{M \in \mathbb{M}_{n,1}} \bC \cdot X^M \right\}
$$
of $\mcal{A}(\mathrm{Mat}_n)$,
and consider the tensor product
$$
\mathrm{ML}_n^{\bullet}\otimes_\C
\kakko{\bigoplus_{k=0}^\infty \bC \cdot w(\alpha,s-k)}
$$
where $\{ w(\alpha,s-k)\}_{k \ge 0}$ are formal vectors.
We introduce a $\U$-module structure on it by
\begin{equation} \label{actionComplexCase}
Y \cdot (F \otimes w(\alpha,s-k)) = 
(Y \cdot F) \otimes w(\alpha,s-k) + (s-k) F (Y \cdot \dete{\alpha}(X) ) \otimes w(\alpha,s-k-1)
\end{equation}
for $Y \in \gl_n$ and $F \in \mathrm{ML}_n^{\bullet}$.
Let $\mcal{ML}_n(\alpha,s)$ be the quotient $\U$-module of 
$\mathrm{ML}_n^{\bullet} \otimes_{\bC} \( \bigoplus_{k=0}^\infty \bC \cdot w(\alpha,s-k) \)$
with respect to the submodule generated by
\begin{equation}\label{eq:relations_ML}
(F \cdot \dete{\alpha}(X)) \otimes w(\alpha,s-k)%
-F \otimes w(\alpha,s-k+1) \quad
(F \in \mathrm{ML}_n^{\bullet},\ w(\alpha,0)=1).
\end{equation}
For $F(X)\in\mathrm{ML}_n^{\bullet}$ and $k\in\Znn$,
we denote by $F(X)\adet(X)^{s-k}$
the element in $\mcal{ML}_n(\alpha,s)$
represented by $F(X)\otimes w(\alpha,s-k)$.
We notice that
$\adet(X)\adet(X)^{s-k}=\adet(X)^{s-k+1}$
by \eqref{eq:relations_ML}.

Denote by $\mcal{V}(\alpha,s)$ the submodule of $\mcal{ML}_n(\alpha,s)$,
generated by the vector $\adet(X)^s(=1\otimes w(\alpha,s))$.
When $s$ is a non-negative integer $l$,
we can naturally consider that  
$$
\mcal{ML}_n(\alpha,l) \subset \mathrm{ML}_n^{\bullet}
$$
because $1 \otimes w(\alpha,l) = \dete{\alpha}(X)^l \otimes 1$
so that it follows that
$$
\mcal{V}(\alpha,l)\cong\U\cdot\dete{\alpha}(X)^l.
$$
Thus we regard the space
$\mcal{V}(\alpha,s)$
as a suitable formulation
of the cyclic module $\U\cdot\dete{\alpha}(X)^s$
for $s\in\C$.
We note that $\mcal{ML}_n(\alpha,l)$ can be realized
in the quotient field for the algebra $\mathrm{ML}_n^{\bullet}$
when $l$ is a negative integer.

\begin{example}
Let $\alpha=-1$.
Then 
$\mcal{V}(-1,s)=\U \cdot \det(X)^s$  is one-dimensional space 
and we have $E_{pp} \cdot\adet[-1](X)^s=s\adet[-1](X)^s$
for any $1 \le p \le n$ from \eqref{actionComplexCase}.
Thus the module $\mcal{V}(-1,s)$ is the irreducible module
with ``highest weight $(s,s,\dots,s)$''. \qed
\end{example}

The module $\mcal{V}(\alpha,s)$ is infinite dimensional in general.
For instance,
if $\alpha=0$ and $s \in \bC\setminus\bZ_{+}$,
then we see that
$$
E_{12}^k \cdot \adet[0](X)^s=s(s-1)\cdots(s-k+1) 
(x_{11} x_{12} x_{33} x_{44} \cdots x_{nn})^k \adet[0](X)^{s-k}
$$
for each $k \ge 0$ and these vectors are linearly independent,
and this obviously implies $\dim_\C\mcal{V}(0,s)=\infty$.
In the infinite dimensional cases,
the following two problems are fundamental:
\begin{enumerate}
\item Unitarizability of each irreducible subrepresentation
appearing in the decomposition of the cyclic module $\U\cdot\adet(X)^s$.
\item Description of the ``\emph{content function}'' for
each isotypic component in $\U\cdot\adet(X)^s$ 
as a certain special function such as a solution of some Fuchsian
type ordinary differential equation. (See Remark \ref{remark:Heun}.)
\end{enumerate}
We will treat these problems in our future studies.

\subsection{Generalized immanants}

Let $\phi$ be a class function on $\mf{S}_n$.
We define the {\it $\phi$-immanant} by
$$
\imm{\phi}(X)= \sum_{\sigma \in \mf{S}_n} \phi(\sigma) x_{1 \sigma(1)} \cdots x_{n \sigma(n)}.
$$
For $l$ class functions $\phi_1,\dots,\phi_l$, consider the cyclic module
$$
\U\cdot\prod_{i=1}^l\imm{\phi_i}(X),
$$
which is the submodule of $\bigoplus_{M \in \mathbb{M}_{n,l}} \bC \cdot X^M$.
In the article we discuss
the special case where $\phi_1(\sigma) = \cdots =\phi_l(\sigma)=\alpha^{\nu(\sigma)}$.
The discussion, and hence several propositions, in \S 3 can be extended to this generalized situation
because we do not use any special feature of the function
$\alpha^{\nu(\sigma)}$. See the Appendix below.

\begin{ackn}
The authors thank Jyoichi Kaneko for fruitful discussion on the Jacobi polynomials.
\end{ackn}


\def\M{\mathbb{M}}
\def\GLmod#1#2{\mathcal{M}_{#1}^{#2}}
\def\labelenumi{{\upshape(\theenumi)}}

\section{Appendix: Transition matrices and zonal spherical functions}
\begin{center}
\scshape by Kazufumi KIMOTO
\end{center}

We investigate the structure of the cyclic module
$\Cmod nl(\alpha)\deq\U\cdot\adet(X)^l$ by embedding it to the tensor product space
$(\C^n)^{\otimes nl}$ and utilizing the Schur-Weyl duality.
We show that the entries of the transition matrices $\tramat\lambda(\alpha)$
are given by a variation of
the spherical Fourier transformation
of a certain class function on $\sym{nl}$
with respect to the subgroup $\sym l^n$ (Theorem \ref{thm:mymain}).
This result also provides another proof of Theorem \ref{thm:matrixF}.
Further, we calculate the polynomial $\tramat[2,l]{(2l-s,s)}(\alpha)$
by using an explicit formula of the values of zonal spherical functions
for the Gelfand pair $(\sym{2n},\sym n\times\sym n)$ due to Bannai and Ito
(Theorem \ref{thm:my_n=2_case}).

\subsection{Irreducible decomposition of $\Cmod nl(\alpha)$ and transition matrices}

Fix $n,l\in\N$.
Consider the standard tableau $\reftab$ with shape $(l^n)$
such that the $(i,j)$-entry of $\reftab$ is $(i-1)l+j$.
For instance, if $n=3$ and $l=2$, then
\begin{equation*}
\reftab=\lower1.5em\hbox{\setlength{\unitlength}{1.15em}%
\begin{picture}(2,3)
\multiput(0,0)(0,1){4}{\line(1,0){2}}
\multiput(0,0)(1,0){3}{\line(0,1){3}}
\put(0,2){\makebox(1,1){$1$}}
\put(1,2){\makebox(1,1){$2$}}
\put(0,1){\makebox(1,1){$3$}}
\put(1,1){\makebox(1,1){$4$}}
\put(0,0){\makebox(1,1){$5$}}
\put(1,0){\makebox(1,1){$6$}}
\end{picture}}\,.
\end{equation*}
We denote by $K=R(\reftab)$ and $H=C(\reftab)$
the row group and column group of the standard tableau $\reftab$
respectively.
Namely,
\begin{equation}
K=\ckakko{g\in\sym{nl}\,\big|\,\ceil{g(x)}=\ceil{x},\ x\in[nl]},\quad
H=\ckakko{g\in\sym{nl}\,\big|\,g(x)\equiv x\pmod l,\ x\in[nl]}.
\end{equation}
We put
\begin{equation}
e=\frac1{\card K}\sum_{k\in K}k\in\C[\sym{nl}].
\end{equation}
This is clearly an idempotent element in $\C[\sym{nl}]$.
Let $\varphi$ be a class function on $H$.
We put
\begin{align*}
\Phi\deq\sum_{h\in H}\varphi(h)h\in\C[\sym{nl}].
\end{align*}
Consider the tensor product space $V=(\C^n)^{\otimes nl}$.
We notice that $V$ has a $(\U,\C[\sym{nl}])$-module structure given by
\begin{align*}
E_{ij}\cdot \ve_{i_1}\otimes\dots\otimes\ve_{i_{nl}}
&\deq
\sum_{s=1}^{nl}\delta_{i_{s},j}\,
\ve_{i_{1}}\otimes\dots\otimes\overset{\text{$s$-th}}{\ve_i}\otimes\dots\otimes\ve_{i_{nl}},\\
\ve_{i_1}\otimes\dots\otimes\ve_{i_{nl}}\cdot\sigma
&\deq
\ve_{i_{\sigma(1)}}\otimes\dots\otimes\ve_{i_{\sigma(nl)}}
\qquad(\sigma\in\sym{nl})
\end{align*}
where $\{\ve_i\}_{i=1}^{n}$ denotes the standard basis of $\C^n$.
The main concern of this subsection is to describe the irreducible decomposition of
the left $\U$-module $V\cdot e\Phi e$.

We first show that $\Cmod nl(\alpha)$ is isomorphic to $V\cdot e\Phi e$
for a special choice of $\varphi$.
Consider the group isomorphism
$\theta:H\to\sym n^l$ defined by
\begin{align*}
\theta(h)\deq(\theta(h)_1,\dots,\theta(h)_l);\quad
\theta(h)_i(x)=y \iff h((x-1)l+i)=(y-1)l+i.
\end{align*}
We also define an element $D(X;\phi)\in\A$ by
\begin{align*}
D(X;\phi)&\deq
\sum_{h\in H}\phi(h)
\prod_{q=1}^n\prod_{p=1}^l x_{\theta(h)_p(q),q}
=\sum_{h\in H}\phi(h)
\prod_{q=1}^n\prod_{p=1}^l x_{q,\theta(h)^{-1}_p(q)}\\
&=\sum_{\sigma_1,\dots,\sigma_l\in\sym n}
\phi(\theta^{-1}(\sigma_1,\dots,\sigma_l))
\prod_{q=1}^n\prod_{p=1}^l x_{\sigma_p(q),q}.
\end{align*}
We note that $D(X;\anu\cdot)=\adet(X)^l$
since
$\cw{\theta^{-1}(\sigma_1,\dots,\sigma_l)}
=\cw{\sigma_1}+\dots+\cw{\sigma_l}$
for $(\sigma_1,\dots,\sigma_l)\in\sym l^n$.

Take a class function $\delta_H$ on $H$ defined by
\begin{align*}
\delta_H(h)=\begin{cases}
1 & h=1 \\ 0 & h\ne1.
\end{cases}
\end{align*}
We see that $D(X;\delta_H)=(x_{11}x_{22}\dots x_{nn})^l$.
We need the following lemma
(The assertion (1) is just a rewrite of Lemma \ref{lem:Cyclic},
and (2) is immediate to verify).

\begin{lem}\label{lem:cyclic}
{\upshape(1)}
It holds that
\begin{align*}
&\U\cdot\ve_1^{\otimes l}\otimes\dots\otimes\ve_n^{\otimes l}
=V\cdot e=\st^l(\C^n)^{\otimes n},\\
&\U\cdot D(X;\delta_H)
=\bigoplus_{\substack{i_{pq}\in\{1,2,\dots,n\}\\(1\le p\le l,\,1\le q\le n)}}
\C\cdot\prod_{q=1}^n\prod_{p=1}^l x_{i_{pq}q}
\cong\st^l(\C^n)^{\otimes n}.
\end{align*}

\noindent{\upshape(2)}
The map
\begin{align*}
\intertwiner:\U\cdot D(X;\delta_H)\ni
\prod_{q=1}^n\prod_{p=1}^l x_{i_{pq}q}
\longmapsto
(\ve_{i_{11}}\otimes\dots\otimes\ve_{i_{l1}})\otimes\dots\otimes%
(\ve_{i_{1n}}\otimes\dots\otimes\ve_{i_{ln}})\cdot e\in V\cdot e
\end{align*}
is a bijective $\U$-intertwiner.
\qed
\end{lem}

We see that
\begin{align*}
\intertwiner\kakko{D(X;\phi)}
&=\sum_{h\in H}\phi(h)
\intertwiner\kakko{\prod_{q=1}^n\prod_{p=1}^l x_{\theta(h)_p(q),q}}\\
&=\sum_{h\in H}\phi(h)
(\ve_{\theta(h)_1(1)}\otimes\dots\otimes\ve_{\theta(h)_l(1)})\otimes\dots\otimes%
(\ve_{\theta(h)_1(n)}\otimes\dots\otimes\ve_{\theta(h)_l(n)})\cdot e\\
&=\ve_1^{\otimes l}\otimes\dots\otimes\ve_n^{\otimes l}%
\cdot\sum_{h\in H}\phi(h)h\cdot e
=\ve_1^{\otimes l}\otimes\dots\otimes\ve_n^{\otimes l}%
\cdot e\Phi e
\end{align*}
by (2) in Lemma \ref{lem:cyclic}.
Using (1) in Lemma \ref{lem:cyclic},
we have the
\begin{lem}
It holds that
\begin{align*}
\U\cdot D(X;\phi)\cong V\cdot e\Phi e
\end{align*}
as a left $\U$-module.
In particular, $V\cdot e\Phi e\cong\Cmod nl(\alpha)$
if $\varphi(h)=\anu h$.
\qed
\end{lem}



By the Schur-Weyl duality, we have
\begin{align*}
V\cong\bigoplus_{\lambda\vdash nl}\GLmod n\lambda\boxtimes\Smod\lambda.
\end{align*}
Here
$\Smod\lambda$ denotes the irreducible unitary right $\sym{nl}$-module
corresponding to $\lambda$.
We see that
\begin{align*}
\dim\kakko{\Smod\lambda\cdot e}
=\inprod[\sym{nl}]{\ind_K^G\triv K}{\Smod\lambda}
=K_{\lambda(l^n)},
\end{align*}
where $\triv K$ is the trivial representation of $K$ and
$\inprod[\sym{nl}]\pi\rho$ is the intertwining number
of given representations $\pi$ and $\rho$ of $\sym{nl}$.
Since $K_{\lambda(l^n)}=0$ unless $\len\lambda\le n$,
it follows the
\begin{thm}
It holds that
\begin{align*}
V\cdot e\Phi e
\cong\bigoplus_{\substack{\lambda\vdash nl\\ \len\lambda\le n}}
\GLmod n\lambda\boxtimes\kakko{\Smod\lambda\cdot e\Phi e}.
\end{align*}
In particular, as a left $\U$-module,
the multiplicity of $\GLmod n\lambda$ in $V\cdot e\Phi e$
is given by
\begin{align*}
\dim\kakko{\Smod\lambda\cdot e\Phi e}
=\rank_{\End(\Smod\lambda\cdot e)}(e\Phi e).
\end{align*}
\qed
\end{thm}

Let $\lambda\vdash nl$ be a partition such that $\len\lambda\le n$
and put $d=K_{\lambda(l^n)}$.
We fix an orthonormal basis
$\{\ve_1^\lambda,\dots,\ve_{f^\lambda}^\lambda\}$
of $\Smod\lambda$ such that the first $d$ vectors
$\ve_1^\lambda,\dots,\ve_d^\lambda$ form a subspace
$(\Smod\lambda)^K$ consisting of $K$-invariant vectors
and left $f^\lambda-d$ vectors form
the orthocomplement of $(\Smod\lambda)^K$
with respect to the $\sym{nl}$-invariant inner product.
The matrix coefficient of $\Smod\lambda$
relative to this basis is
\begin{equation}\label{eq:def_of_psi}
\psi^\lambda_{ij}(g)
=\inprod[\Smod\lambda]{\ve_i^\lambda\cdot g}{\ve_j^\lambda}
\quad(g\in\sym{nl},\ 1\le i,j\le f^\lambda).
\end{equation}
We notice that this function is $K$-biinvariant.
We see that the multiplicity of $\GLmod n\lambda$ in $V\cdot e\Phi e$
is given by the rank of the matrix
\begin{align*}
\kakko{\sum_{h\in H}\varphi(h)\psi^\lambda_{ij}(h)}_{1\le i,j\le d}.
\end{align*}

As a particular case, we obtain the
\begin{thm}\label{thm:mymain}
The multiplicity of the irreducible representation $\GLmod n\lambda$
in the cyclic module $\U\cdot\adet(X)^l$ is equal to the rank of
\begin{equation}
\tramat\lambda(\alpha)
=\kakko{\sum_{h\in H}\anu h\psi^\lambda_{ij}(h)}_{1\le i,j\le d},
\end{equation}
where $\{\psi^\lambda_{ij}\}_{i,j}$ denotes a basis of
the $\lambda$-component of the space $C(K\backslash\sym{nl}/K)$
of $K$-biinvariant functions on $\sym{nl}$
given by \eqref{eq:def_of_psi}.
\end{thm}

\begin{rem}
\begin{enumerate}
\item
By the definition of the basis $\{\psi^\lambda_{ij}\}_{i,j}$
in \eqref{eq:def_of_psi},
we have $\tramat\lambda(0)=I$.
\item
Since $\anu{g^{-1}}=\anu g$ and
$\psi^\lambda_{ij}(g^{-1})=\overline{\psi^\lambda_{ji}(g)}$
for any $g\in\sym{nl}$,
the transition matrices satisfy
$\tramat\lambda(\alpha)^*=\tramat\lambda(\overline\alpha)$.
\item
In Examples \ref{ex:MW2005} and \ref{ex:standard_case_again} below,
the transition matrices
are given by \emph{diagonal matrices}.
We expect that any transition matrix $\tramat\lambda(\alpha)$
is \emph{diagonalizable} in $\Mat_{K_{\lambda(l^n)}}(\C[\alpha])$.
\end{enumerate}
\end{rem}

\begin{ex}\label{ex:MW2005}
If $l=1$, then $H=G=\sym n$ and $K=\{1\}$.
Therefore, for any $\lambda\vdash n$, we have
\begin{equation}
\tramat[n,1]\lambda(\phi)
=\frac{n!}{f^\lambda}\inprod[\sym n]{\phi}{\chi^\lambda}I
\end{equation}
by the orthogonality of the matrix coefficients. 
Here $\chi^\lambda$ denotes the irreducible character of $\sym n$
corresponding to $\lambda$.
In particular, if $\phi=\anu\cdot$, then
\begin{equation}
\tramat[n,1]\lambda(\alpha)
=\cp\lambda(\alpha)I
\end{equation}
since the Fourier expansion of $\anu\cdot$ 
(as a class function on $\sym n$) is
\begin{equation}\label{eq:specialFCF}
\anu\cdot
=\sum_{\lambda\vdash n}
\frac{f^\lambda}{n!}f_\lambda(\alpha)\chi^\lambda,
\end{equation}
which is obtained by specializing
the Frobenius character formula for $\sym n$
(see, e.g. \cite{Mac}).
\end{ex}

\begin{ex}
Let us calculate $\tramat{(nl)}(\alpha)$ by using Theorem \ref{thm:mymain}.
Since $\Smod{(nl)}$ is the trivial representation, it follows that
$(\Smod{(nl)})^K=\Smod{(nl)}$ and
\begin{align*}
\tramat{(nl)}(\alpha)=
\sum_{h\in H}\anu h\inprod{\ve\cdot h}{\ve}=\sum_{\sigma_1,\dots,\sigma_l\in\sym n}\anu{\sigma_1}\dots\anu{\sigma_l}
=\kakko{(1+\alpha)(1+2\alpha)\dots(1+(n-1)\alpha)}^l,
\end{align*}
where $\ve$ denotes a unit vector in $\Smod{(nl)}$.
\end{ex}

\begin{ex}\label{ex:standard_case_again}
\def\stirling#1#2{\genfrac{[}{]}{0pt}{}{\,#1\,}{\,#2\,}}
Let us calculate $\tramat{(nl-1,1)}(\alpha)$ by using Theorem \ref{thm:mymain}.
As is well known,
the irreducible (right) $\sym{nl}$-module $\Smod{(nl-1,1)}$ can be realized in $\C^{nl}$ as follows:
\begin{align*}
\Smod{(nl-1,1)}=\Biggl\{(x_j)_{j=1}^{nl}\in\C^{nl}\ \Bigg|\ \sum_{j=1}^{nl}x_j=0\Biggr\}.
\end{align*}
This is a unitary representation with respect to the ordinary hermitian inner product $\inprod\cdot\cdot$
on $\C^{nl}$.
It is immediate to see that
\begin{align*}
\kakko{\Smod{(nl-1,1)}}^{\!K}=\ckakko{(x_j)_{j=1}^{nl}\in\Smod{(nl-1,1)}\,\Big|\,
x_{pl+1}=x_{pl+2}=\dots=x_{(p+1)l}\quad(0\le p<n)}.
\end{align*}
Take an orthonormal basis $\ve_1,\dots,\ve_{n-1}$ of $\kakko{\Smod{(nl-1,1)}}^K$ by
\begin{align*}
\ve_j=\frac1{\sqrt{nl}}
\bigl(\overbrace{\omega^j,\dots,\omega^j}^l,
\overbrace{\omega^{2j},\dots,\omega^{2j}}^l,
\dots,
\overbrace{\omega^{nj},\dots,\omega^{nj}}^l\bigr)
\qquad(1\le j\le n-1),
\end{align*}
where $\omega$ is a primitive $n$-th root of unity.
Then, the $(i,j)$-entry of the transition matrix $\tramat{(nl-1,1)}(\alpha)$ is
\begin{align*}
\sum_{h\in H}\anu h\inprod{\ve_i\cdot h}{\ve_j}
&=\frac1{nl}\sum_{\sigma_1,\dots,\sigma_l\in\sym n}\sum_{p=1}^n\sum_{q=1}^l
\anu{\sigma_1}\dots\anu{\sigma_l}\omega^{\sigma_q(p)i-pj}\\
&=\kakko{\sum_{\tau\in\sym n}\anu\tau}^{l-1}
\kakko{\frac1n\sum_{\sigma\in\sym n}\sum_{p=1}^n\anu\sigma\omega^{\sigma(p)i-pj}}.
\end{align*}
The first factor is $\kakko{(1+\alpha)(1+2\alpha)\dots(1+(n-1)\alpha)}^{l-1}$.
We show that
\begin{align*}
\frac1n\sum_{\sigma\in\sym n}\sum_{p=1}^n\anu\sigma\omega^{\sigma(p)i-pj}
=(1-\alpha)(1+\alpha)(1+2\alpha)\dots(1+(n-2)\alpha)\delta_{ij}
\qquad(i,j=1,2,\dots,n-1).
\end{align*}
For this purpose,
by comparing the coefficients of $\alpha^{n-m}$ in both sides,
it is enough to prove
\begin{align*}
\frac1n\sum_{\substack{\sigma\in\sym n\\\nu(\sigma)=n-m}}\sum_{p=1}^n\omega^{\sigma(p)i-pj}
=\ckakko{\stirling{n-1}{m-1}-\stirling{n-1}m}\delta_{ij}
\qquad(i,j,m=1,2,\dots,n-1),
\end{align*}
where $\stirling nm$ denotes the Stirling number of the first kind
(see, e.g. \cite{GKP} for the definition).
Since
\begin{align*}
\#\ckakko{\sigma\in\sym n\,;\,\nu(\sigma)=n-m,\ \sigma(p)=x}=
\begin{cases}
\stirling{n-1}{m-1} & x=p,\\
\stirling{n-1}{m} & x\ne p
\end{cases}
\end{align*}
for each $p,x\in[n]$, it follows that
\begin{align*}
\frac1n\sum_{\sigma\in\sym n}\sum_{p=1}^n\anu\sigma\omega^{\sigma(p)i-pj}
&=\frac1n\sum_{p=1}^n\omega^{-pj}
\ckakko{\stirling{n-1}{m-1}\omega^{pi}+\sum_{x\ne p}\stirling{n-1}m\omega^{xi}}\\
&=\ckakko{\stirling{n-1}{m-1}-\stirling{n-1}m}\frac1n\sum_{p=1}^n\omega^{p(i-j)}\\
&=\ckakko{\stirling{n-1}{m-1}-\stirling{n-1}m}\delta_{ij},
\end{align*}
which is the required conclusion.
Here we notice that $\sum_{x\ne p}\omega^{xi}=-\omega^{pi}$ since $1\le i<n$.
Consequently, we obtain
\begin{align*}
\tramat{(nl-1,1)}(\alpha)
=\kakko{(1-\alpha)\kakko{(1+\alpha)(1+2\alpha)\dots(1+(n-2)\alpha)}^{l}(1+(n-1)\alpha)^{l-1}\delta_{ij}}_{1\le i,j\le n-1},
\end{align*}
so that the multiplicity of $\GLmod n{(nl-1,1)}$ in $\Cmod nl(\alpha)$ is
zero if $\alpha\in\{1,-1,-1/2,\dots,-1/(n-1)\}$ and
$n-1$ otherwise.
\end{ex}

The trace of the transition matrix $\tramat\lambda(\alpha)$ is
\begin{equation}
\gcp\lambda{n,l}\deq
\tr\tramat\lambda(\alpha)
=\sum_{h\in H}\anu h\omega^\lambda(h),
\end{equation}
where $\omega^\lambda$ is the \emph{zonal spherical function} for $\lambda$
with respect to $K$ defined by
\begin{align*}
\omega^\lambda(g)\deq
\frac1{\card K}\sum_{k\in K}\chi^\lambda(kg)
\quad(g\in\sym{nl}).
\end{align*}
This is regarded as a generalization of the modified content polynomial
since $\gcp\lambda{n,1}=f^\lambda f_\lambda(\alpha)$ as we see above.
It is much easier to handle these polynomials
than the transition matrices.
If we could prove that a transition matrix $\tramat\lambda(\alpha)$ is a scalar matrix,
then we would have $\tramat\lambda(\alpha)=d^{-1}\gcp\lambda{n,l}I$
($d=\dim(\Smod\lambda)^K$)
and hence we see that the multiplicity of $\GLmod n\lambda$
in $\Cmod nl(\alpha)$ is completely controlled
by the single polynomial $\gcp\lambda{n,l}$.
In this sense,
it is desirable to obtain a characterization
of the irreducible representations
whose corresponding transition matrices are scalar
as well as to get an explicit expression for the polynomials $\gcp\lambda{n,l}$.
Here we give a sufficient condition for $\lambda\vdash nl$ such that
$\tramat\lambda(\alpha)$ is a scalar matrix.

\begin{prop}
\begin{enumerate}
\item Denote by $N_H(K)$ the normalizer of $K$ in $H$.
The transition matrix $\tramat\lambda(\alpha)$ is scalar
if $(\Smod\lambda)^K$ is irreducible as a $N_H(K)$-module.
\item
If $\lambda$ is of hook-type $($i.e. $\lambda=(nl-r,1^r)$ for some $r<n$$)$,
then $\tramat\lambda(\alpha)$ is scalar.
\end{enumerate}
\end{prop}

\begin{proof}
\def\vx{\boldsymbol{x}}
Notice that $N_H(K)\cong\sym n$.
Consider a linear map $T\in\End((\Smod\lambda)^K)$ given by
\begin{align*}
T(\vx)=\sum_{j=1}^d\kakko{\sum_{h\in H}\anu h\inprod[\Smod\lambda]{\vx\cdot h}{\ve^\lambda_j}}\ve^\lambda_j
\qquad(\vx\in(\Smod\lambda)^K),
\end{align*}
where $d=\dim(\Smod\lambda)^K$.
It is direct to check that $T$ gives an intertwiner of $(\Smod\lambda)^K$ as a $N_H(K)$-module.
Hence, by Schur's lemma, $T$ is a scalar map (and $\tramat\lambda(\alpha)$ is a scalar matrix)
if $(\Smod\lambda)^K$ is an irreducible $N_H(K)$-module.
When $\lambda=(nl-r,1^r)$ for some $r<n$,
it is proved in \cite[Proposition 5.3]{AMT} that
$(\Smod{(nl-r,1^r)})^K\cong\Smod{(n-r,1^r)}$ as $N_H(K)$-modules.
Thus we have the proposition.
\end{proof}

\begin{ex}
Let us calculate $\gcp{(nl-1,1)}{n,l}$.
We notice that $\chi^{(nl-1,1)}(g)=\fix_{nl}(g)-1$
where $\fix_{nl}$ denotes the number of fixed points
in the natural action $\sym{nl}\curvearrowright[nl]$.
Hence we see that
\begin{align*}
\gcp{(nl-1,1)}{n,l}
&=\sum_{h\in H}\anu h\frac1{\card K}\sum_{k\in K}(\fix_{nl}(kh)-1)
=\sum_{h\in H}\anu h\frac1{\card K}\sum_{k\in K}\sum_{x\in[nl]}\delta_{khx,x}%
-\sum_{h\in H}\anu h.
\end{align*}
It is easily seen that $khx\ne x$ for any $k\in K$ if $hx\ne x$ ($x\in[nl]$).
Thus it follows that
\begin{align*}
\frac1{\card K}\sum_{k\in K}\sum_{x\in[nl]}\delta_{khx,x}
=\sum_{x\in[nl]}\delta_{hx,x}\frac1{\card K}\sum_{k\in K}\delta_{kx,x}
=\frac1l\fix_{nl}(h)\qquad(h\in H).
\end{align*}
Therefore we have
\begin{align*}
\gcp{(nl-1,1)}{n,l}
&=\frac1l\sum_{h\in H}\anu h\fix_{nl}(h)-\sum_{h\in H}\anu h
=\gcp{(n)}{n,1}^{l-1}\gcp{(n-1,1)}{n,1}\\
&=(n-1)(1-\alpha)(1-(n-1)\alpha)^{l-1}\prod_{i=1}^{n-2}(1+i\alpha)^l.
\end{align*}
Since the transition matrix $\tramat{(nl-1,1)}$ is a scalar one and
its size is $\dim\Smod{(n-1,1)}=n-1$, we get
$\tramat{(nl-1,1)}(\alpha)=(1-\alpha)(1-(n-1)\alpha)^{l-1}\prod_{i=1}^{n-2}(1+i\alpha)^lI_{n-1}$ again.
\end{ex}

We will investigate these polynomials $\gcp\lambda{n,l}$
and their generalizations in \cite{K2007c}.

\subsection{Irreducible decomposition of $\Cmod 2l(\alpha)$ and Jacobi polynomials}

In this section,
as a particular example,
we consider the case where $n=2$ and
calculate the transition matrix $\tramat[2,l]\lambda(\alpha)$ explicitly.
Since the pair $(\sym{2l},K)$ is a \emph{Gelfand pair}
(see, e.g. \cite{Mac}),
it follows that
\begin{align*}
K_{\lambda(l^2)}
=\inprod[\sym{2l}]{\ind_K^{\sym{2l}}\triv K}{\Smod\lambda}
=1
\end{align*}
for each $\lambda\vdash2n$ with $\len\lambda\le2$.
Thus, in this case,
the transition matrix is just a polynomial and is given by
\begin{equation}
\tramat[2,l]\lambda(\alpha)
=\tr\tramat[2,l]\lambda(\alpha)
=\sum_{h\in H}\anu h\omega^\lambda(h)
=\sum_{s=0}^l\binom ls\omega^\lambda(g_s)\alpha^s.
\end{equation}
Here we put $g_s=(1,l+1)(2,l+2)\dots(s,l+s)\in\sym{2n}$.
Now we write $\lambda=(2l-p,p)$ for some $p$ ($0\le p\le l$).
The value $\omega^{(2l-p,p)}(g_s)$ of the zonal spherical function is
calculated by Bannai and Ito \cite[p.218]{BI} as
\begin{align*}
\omega^{(2l-p,p)}(g_s)
=Q_p(s;-l-1,-l-1,l)
=\sum_{j=0}^p(-1)^j\binom pj\binom{2l-p+1}j\binom lj^{\!\!-2}
\binom sj,
\end{align*}
where
\begin{align*}
Q_n(x;\alpha,\beta,N)&\deq
\tHGF32{-n,n+\alpha+\beta+1,-x}{\alpha+1,-N}{1}\\
&=\sum_{j=0}^N(-1)^j
\binom nj\binom{-n-\alpha-\beta-1}j
\binom{-\alpha-1}j^{\!\!-1}\binom Nj^{\!\!-1}
\binom xj
\end{align*}
is the Hahn polynomial
(see also \cite[p.399]{Mac}).
We also denote by $\tHGF{n+1}{n}{a_1,\dots,a_{p}}{b_1,\dots,b_{q-1},-N}x$
the hypergeometric polynomial
\begin{align*}
\tHGF pq{a_1,\dots,a_p}{b_1,\dots,b_{q-1},-N}x
=\sum_{j=0}^N\frac{\phs{a_1}j\dots\phs{a_{p}}j}
{\phs{b_1}j\dots\phs{b_{q-1}}j\phs{-N}j}\frac{x^j}{j!}
\end{align*}
for $p,q,N\in\N$ in general (see \cite{AAR}).
We now re-prove Theorem \ref{thm:CaseTwo} as follows:
\begin{thm}\label{thm:my_n=2_case}
Let $l$ be a positive integer.
It holds that
\begin{align*}
\tramat[2,l]{(2l-p,p)}(\alpha)
=\sum_{s=0}^l\binom ls Q_p(s;l-1,l-1,l)\alpha^s
=(1+\alpha)^{l-p}G_p^l(\alpha)
\end{align*}
for $p=0,1,\dots,l$.
\end{thm}

\begin{proof}
Let us put $x=-1/\alpha$.
Then we have
\begin{align*}
\sum_{s=0}^l\binom ls Q_p(s;l-1,l-1,l)\alpha^s
&=\sum_{j=0}^p(-1)^j\binom pj\binom{2l-p+1}j\binom lj^{\!\!-1}
\alpha^j(1+\alpha)^{l-j}\\
&=x^{-l}(x-1)^{l-p}\sum_{j=0}^p\binom pj\binom{2l-p+1}j\binom lj^{\!\!-1}(x-1)^{p-j}
\end{align*}
and
\begin{align*}
(1+\alpha)^{l-p}G_p^l(\alpha)
&=x^{-l}(x-1)^{l-p}\sum_{j=0}^p(-1)^j
\binom pj\binom{l-p+j}j\binom lj^{\!\!-1}(-x)^{p-j}.
\end{align*}
Here we use the elementary identity
\begin{equation*}
\sum_{s=0}^l \binom ls\binom sj\alpha^s=\binom lj\alpha^j(1+\alpha)^{l-j}.
\end{equation*}
Hence, to prove the theorem,
it is enough to verify
\begin{equation}
\sum_{i=0}^p\binom pi\binom{l-p+i}i\binom li^{\!\!-1}x^{p-i}
=\sum_{j=0}^p\binom pj\binom{2l-p+1}j\binom lj^{\!\!-1}(x-1)^{p-j}.
\end{equation}
Comparing the coefficients of Taylor expansion
of these polynomials at $x=1$,
we notice that the proof is reduced to the equality
\begin{equation}
%
\sum_{i=0}^r\binom{l-i}{l-r}\binom{l-p+i}{l-p}
=\binom{2l-p+1}r
\end{equation}
for $0\le r\le p$, which is well known (see, e.g. (5.26) in \cite{GKP}).
Hence we have the conclusion.
\end{proof}
Thus we obtain the irreducible decomposition
\begin{align}
\Cmod 2l(-1)\cong\GLmod2{(l,l)},\qquad
\Cmod 2l(\alpha)\cong\bigoplus_{\substack{0\le p\le l\\G_p^l(\alpha)\ne0}}\GLmod2{(2l-p,p)}
\quad(\alpha\ne-1)
\end{align}
of $\Cmod2l(\alpha)$ again.

\begin{rem}
\begin{enumerate}
\item The calculation above uses the advantage for the fact
that $(\sym{nl},\sym l^n)$ is the Gelfand pair \emph{only when} $n=2$.
\item We have used the result in \cite[p.218]{BI} for the theorem.
It is worth mentioning that one may prove conversely the result
in \cite[p.218]{BI} from Theorem \ref{thm:CaseTwo}.
\end{enumerate}
\end{rem}

\begin{ackn}
The author would thank
Professor Itaru Terada for noticing that his work \cite{AMT}
is useful for the discussion in \S 6.2.
\end{ackn}


\medskip

\noindent
\textsc{Kazufumi KIMOTO}\\
Department of Mathematical Science, University of the Ryukyus.\\
Nishihara, Okinawa 903-0213, JAPAN.\\
\texttt{kimoto@math.u-ryukyu.ac.jp}\\

\noindent
\textsc{Sho MATSUMOTO}\\
Faculty of Mathematics, Kyushu University.\\
Hakozaki Higashi-ku, Fukuoka, 812-8581 JAPAN.\\

\noindent
\textit{Current address}:\\
Graduate School of Mathematics, Nagoya University.\\
Chikusa, Nagoya 464-8602, JAPAN.\\
\texttt{sho-matsumoto@math.nagoya-u.ac.jp}\\

\noindent
\textsc{Masato WAKAYAMA}\\
Faculty of Mathematics, Kyushu University.\\
Hakozaki Higashi-ku, Fukuoka 812-8581, JAPAN.\\
\texttt{wakayama@math.kyushu-u.ac.jp}


\begin{thebibliography}{AAR}
%
%
\bibitem[AAR]{AAR}
G.\,E. Andrews, R. Askey and R. Roy,
``Special Functions'',
Encyclopedia of Mathematics and its Applications, 71.
Cambridge University Press, Cambridge, 1999.

\bibitem[AMT]{AMT}
S. Ariki, J. Matsuzawa and I. Terada,
{\it Representation of Weyl groups on zero weight spaces of $\mathfrak{g}$-modules.}
Algebraic and topological theories (Kinosaki, 1984),
546--568, Kinokuniya, Tokyo, 1986.

\bibitem[BI]{BI}
E. Bannai and T. Ito,
``Algebraic Combinatorics I, Association Schemes.''
The Benjamin/Cummings Publishing Co., Inc., Menlo Park, CA, 1984.

\bibitem[FH]{FultonHarris}
 W. Fulton and J. Harris, 
 ``Representation theory. A first course'',
  Graduate Texts in Mathematics, 129. Readings in Mathematics. Springer-Verlag, New York, 1991.

\bibitem[GKP]{GKP}
R.\,L. Graham, D.\,E. Knuth and O. Patashnik,
``Concrete Mathematics.
A foundation for computer science'', Second edition.
Addison-Wesley Publishing Company, Reading, MA, 1994.

\bibitem[GN]{GN}
I. M. Gel'fand and M. A. Na{\u\i}mark,
Unitary representations of the Lorentz group,
Acad. Sci. USSR. J. Phys. {\bf10} (1946), 93--94;
Izvestiya Akad. Nauk SSSR. Ser. Mat. {\bf11} (1947), 411--504.

\bibitem[K]{K2007c}
K. Kimoto,
{\it Generalized content polynomials toward $\alpha$-determinant cyclic modules.}
Preprint (2007).

\bibitem[KW1]{KWinvariant}
 K. Kimoto and M. Wakayama,
 {\it Invariant theory for singular $\alpha$-determinants},
 J. Combin. Theory Ser. A {\bf 115} (2008), no. 1, 1--31.

\bibitem[KW2]{KWquantum}
 ------,
 {\it Quantum $\alpha$-determinant cyclic modules of $\mcal{U}_q(\gl_n)$},
 J. Algebra {\bf 313} (2007), 922--956.

\bibitem[Mac]{Mac}
 I. G. Macdonald,
 ``Symmetric Functions and Hall Polynomials'',
 2nd edn., Oxford University Press, 1995.

\bibitem[Mat]{Mat}
 S. Matsumoto, 
 {\it Alpha-pfaffian, pfaffian point process and shifted Schur measure},
 Linear Algebra Appl. {\bf 403} (2005), 369--398.

\bibitem[MW]{MW}
 S. Matsumoto and M. Wakayama,
 {\it Alpha-determinant cyclic modules of $\mathfrak{gl}_n (\mathbb{C})$},
 J. Lie Theory {\bf 16} (2006), 393--405.

\bibitem[S]{Szego}
 G. Szeg\"{o},
 ``Orthogonal Polynomials'', 4th edn., 
  American Mathematical Society, 1975.

\bibitem[SL]{SL}
 S. Yu Slavyanov and W. Lay,
``Special Functions -- A Unified Theory Based on Singularities",
Oxford: Oxford Univ. Press, 2000.

\bibitem[ST]{ST}
 T. Shirai and Y. Takahashi,
 {\it Random point fields associated with certain Fredholm determinants I:
  fermion, Poisson and boson point processes},
 J. Funct. Anal. {\bf 205} (2003), 414--463.

\bibitem[V]{VereJones}
 D. Vere-Jones, 
 {\it A generalization of permanents and determinants},
 Linear Algebra Appl. {\bf 111} (1988), 119--124.

\bibitem[W]{Weyl}
 H. Weyl, 
 ``The Classical Groups. Their invariants and representations'',
  Fifteenth printing. 
  Princeton Landmarks in Mathematics, Princeton University Press, 1997.
\end{thebibliography}
\end{document}